\documentclass[a4paper,11pt]{article}

\usepackage{amsmath}
\usepackage{amsfonts}
\usepackage{amssymb}
\usepackage{amsthm}
\usepackage{graphicx}
\usepackage{psfrag}
\usepackage{subfigure}
\newcommand{\p}{\varphi}
\newcommand{\R}{\mathbb{R}}
\newcommand{\weight}[1]{\langle #1\rangle}

\newcommand{\ZuWeis}{\mathrel{\mathop:\!\!=}}

\newcommand{\D}{\mathcal{D}}
\newcommand{\La}{\mathcal{L}}
\newcommand{\F}{\mathcal{F}}
\newcommand{\K}{\mathcal{K}}
\newcommand{\bv}{\mathbf{v}}
\newcommand{\bw}{\mathbf{w}}
\newcommand{\boldf}{\mathbf{f}}
\newcommand{\tbj}{\widetilde{\mathbf{J}}}
\newcommand{\bpsi}{\boldsymbol{\psi}}

\newcommand{\loc}{\operatorname{loc}}
\newcommand{\uloc}{\operatorname{uloc}}
\newcommand{\N}{\ensuremath{\mathbb{N}}}
\newcommand{\Hone}{H^1_{(0)}}
\newcommand{\supp}{\operatorname{supp}}
\newcommand{\ol}[1]{\overline{#1}}
\newcommand{\sd}{\, d}

\newtheorem{definition}{Definition}[section]
\newtheorem{remark}[definition]{Remark}
\newtheorem{lemma}[definition]{Lemma}
\newtheorem{theorem}[definition]{Theorem}
\newtheorem{assumption}[definition]{Assumption}

\numberwithin{equation}{section}  %equation-counter is set back to zero at the beginning of each section

\renewcommand{\theequation}{\arabic{section}.\arabic{equation}}

\begin{document}
 %%% Neu:
\begin{titlepage}
\title{Existence of Weak Solutions for a Diffuse Interface Model for Two-Phase Flows of 
Incompressible Fluids with Different Densities}
\author{  Helmut Abels\footnote{Fakult\"at f\"ur Mathematik,  
Universit\"at Regensburg,
93040 Regensburg,
Germany, e-mail: {\sf helmut.abels@mathematik.uni-regensburg.de}}, Daniel Depner\footnote{Fakult\"at f\"ur Mathematik,  
Universit\"at Regensburg,
93040 Regensburg,
Germany, e-mail: {\sf daniel.depner@mathematik.uni-regensburg.de}}, and
Harald Garcke\footnote{Fakult\"at f\"ur Mathematik,  
Universit\"at Regensburg,
93040 Regensburg,
Germany, e-mail: {\sf harald.garcke@mathematik.uni-regensburg.de}}}
%\date{}
\end{titlepage}
\maketitle
\begin{abstract}
  We prove existence of weak solutions for a diffuse interface model for the flow of two viscous
 incompressible Newtonian fluids in a bounded domain in two and three space dimensions. In contrast to previous
 works, we study a new model recently developed by Abels, Garcke, and Gr\"un for fluids with different densities, which leads 
 to a solenoidal velocity field. The model is given by a non-homogeneous Navier-Stokes system with a modified convective term coupled to a Cahn-Hilliard
 system. The density of the mixture depends on an order parameter.
\end{abstract}

\noindent{\bf AMS-Classification:} 
Primary: 76T99; %% Two-Phase flows: Others
Secondary: 35Q30, %% Stokes and Navier-Stokes eq.
35Q35, %% Other equations arising in fluid mechanics
76D03, %% Incompressible viscous fluids: Existence, uniqueness, and regularity theory 
76D05, %% Incompressible viscous fluids: Navier-Stokes equations
76D27, %% Incompressible viscous fluids: Other free-boundary flows; Hele-Shaw flows
76D45 %% Incompressible viscous fluids: Capillarity (surface tension)

%%%%%%%%%%%%%%%%%%%%%%%%%%%%%%%%%%%%%%%%%%%%%%%%%%%%%%%%%%%%%%%%%%%%%%%%%%%%%%%%%%%%%%%%%%%%%%%%
%%%%%%%%%%%%%%%%%%%%%%%%%%%%%%%%%%%%%%%%%%%%%%%%%%%%%%%%%%%%%%%%%%%%%%%%%%%%%%%%%%%%%%%%%%%%%%%%
%%%%%%%%%%%%%%%%%%%%%%%%%%%%%%%%%%%%%%%%%%%%%%%%%%%%%%%%%%%%%%%%%%%%%%%%%%%%%%%%%%%%%%%%%%%%%%%%
\section{Introduction} \label{intro}

A fundamental problem in fluid dynamics involves changes in topology
of interfaces between immiscible or partially miscible
fluids. Topological transitions such as pinch off and reconnection of
fluid interfaces are important features of many systems and strongly
affect the flow. Classical models based on sharp interface approaches
typically fail to describe these phenomena. In recent years, diffuse
interface models turned
out to be a promising approach to describe such phenomena.
 In this approach, an order parameter is introduced which
allows for a mixing in an interfacial zone. Therefore  the sharp
interface is replaced by a thin interfacial layer in which the order
parameter, which can be a concentration, rapidly changes its value.
Often these diffuse
interface models (also called {\em phase field models}) allow for
local entropy or free energy inequalities and in such situations they can be
called
thermodynamically consistent. It can be
justified in some cases that sharp interface models are recovered in the
limit
when
the interfacial thickness goes to zero.
We refer to Lowengrub and
Truskinovsky \cite{LT98},
Abels, Garcke and Gr\"un \cite{AGG11} for
a result using formally matched asymptotic expansions and Abels and
R\"oger \cite{AR09}
for a first analytic result. The diffuse interface approach is hence an
attractive approach to model and to numerically simulate fluid
interfaces.

In the literature, diffuse interface models are well established for
two-phase
flows of liquids with identical (``matched'') densities, see Hohenberg and 
Halperin or Gurtin et al.~\cite{HH77, GPV96}.
In contrast, when the
densities are different, several approaches have been discussed in the
literature.
Lowengrub and Truskinovsky \cite{LT98} derived
quasi-incompressible models, where the corresponding velocity field is
not
divergence free. 
On the other hand, Ding et al. \cite{DSS07} proposed a
model with solenoidal fluid velocities which is not known to be thermodynamically
consistent.
Only recently Abels, Garcke and Gr\"un, see \cite{AGG11} derived a thermodynamically
consistent diffuse interface model  for two phase flow
with different densities.
It is the goal of this paper to show existence of weak
solutions for this new model.

More precisely, we consider the following system of Navier-Stokes/Cahn-Hilliard type:
\begin{align*}
 \partial_t (\rho \mathbf{v}) + \operatorname{div} ( \bv \otimes(\rho \bv + \tbj)) - \operatorname{div} (2 \eta(\varphi) D \bv)
  + \nabla p  
  & = - \operatorname{div}(a(\varphi) \nabla \varphi \otimes \nabla \varphi)& \mbox{in } \, Q ,  %\label{newproblem1'} 
\\
 \operatorname{div} \, \bv &= 0& \mbox{in } \, Q,  %\label{newproblem2'} 
\\
 \partial_t \varphi + \bv \cdot \nabla \varphi &= \mbox{div}\left(m(\varphi) \nabla \mu \right)& \mbox{in } \, Q, %\label{newproblem3'} 
\\
 \mu = \Psi'(\varphi) + a'(\varphi) \frac{|\nabla \varphi|^2}{2} &- \operatorname{div}\left( a(\varphi) \nabla \varphi \right)& \mbox{in } \, Q , %\label{newproblem4'} 
\end{align*}
where $\tbj = -\frac{\tilde{\rho}_2 - \tilde{\rho}_1}{2} m(\varphi) \nabla \mu$, 
$Q=\Omega\times(0,\infty)$, and $\Omega \subset \mathbb{R}^d$, $d=2,3$, is a sufficiently smooth bounded domain. We close the system  with the boundary and initial conditions
\begin{alignat*}{2}
 \bv|_{\partial \Omega} &= 0 &\qquad& \text{on}\ \partial\Omega\times (0,\infty), \\
 \partial_n \varphi|_{\partial \Omega} = \partial_n \mu|_{\partial \Omega} &= 0&& \text{on}\ \partial\Omega\times (0,\infty),  \\
 \left(\bv , \varphi \right)|_{t=0} &= \left( \bv_0 , \varphi_0 \right) &&\text{in}\ \Omega, 
\end{alignat*}
where $\partial_n = n\cdot \nabla$ and $n$ denotes the exterior normal at $\partial\Omega$.
Here $\mathbf{v}$ and $\rho=\rho(\varphi)$ are the (mean) velocity  and the
density of the mixture of the two fluids, $p$ is the pressure,
 $\varphi$ is an order parameter related to the concentrations of the two fluids, and $\mu$ is the chemical
potential associated to $\varphi$. 
Moreover,  ${D}\mathbf{v}= \frac12(\nabla \mathbf{v} + \nabla \mathbf{v}^T)$,
$\eta(\varphi)>0$ is a viscosity coefficient, and $m(\varphi)>0$ is a (non-degenerate) mobility coefficient. Furthermore,
$\Psi(\varphi)$ is the homogeneous free energy density for the mixture and the (total) free energy of the system is given by
\begin{equation*}
  \int_{\Omega} \left(a(\varphi)\frac{|\nabla \varphi|^2}2+\Psi(\varphi)\right)\sd x,
\end{equation*}
for some positive coefficient $a(\varphi)$.

An important aspect of our contribution is that we consider a class of singular free energies, 
which includes the homogeneous free energy of the so-called regular solution models used by Cahn and Hilliard~\cite{CahnHilliard}:
\begin{equation}
\label{logpot}\Psi(\varphi) = \frac{\theta}2 \left((1+\varphi)\ln
(1+\varphi)+ (1-\varphi)\ln (1-\varphi)\right) - \frac{\theta_c}2 \varphi^2,\quad \varphi \in [-1,1]
\end{equation}
where
$0<\theta<\theta_c$. Mathematically, these singular free energies ensure that the order parameter stays in the physically 
reasonable interval, which is $[-1,1]$ if $\varphi$ is the difference of volume fractions of both fluids. 
But this leads to singular terms in the equation for the chemical potential. In order to deal with these terms we apply 
techniques, which were developed in Abels and Wilke~\cite{AW07} and applied in Abels~\cite{Abe09b} to the model for matched densities.

In the  case of matched densities, i.\,e., $\tilde{\rho}_1=\tilde{\rho}_2$, $\tbj\equiv 0$, the model reduces 
to the well-known model of matched densities discussed in Hohenberg and Halperin \cite{HH77}. 
In this case existence of weak solutions and well-posedness 
were obtained by Starovoitov~\cite{Sta97}, Boyer~\cite{Boy99}, Liu and Shen~\cite{LS03}, and Abels~\cite{Abe09b}. 
Moreover, Boyer \cite{Boy01} considered a diffuse interface model for fluids with
non-matched densities, which is different from our model. He proved existence of strong solutions, locally in
time, and existence of global weak solutions if the densities of the fluids
are sufficiently close. First analytic results for the model by Lowengrub and Truskinovsky were obtained in \cite{Abe09a,LTModelShortTime}.
For more information on diffuse interface models for two phase flows of incompressible fluids we 
refer to Abels, Garcke and Gr\"un \cite{AGG11}.

The structure of the article is as follows: In Section~\ref{prelimi} we summarize some notation and preliminary results. 
Then, in Section~\ref{secexistence}, we reformulate our system suitably, define weak solutions and state our main result 
on existence of weak solutions. In Section~\ref{sec:Implicit} we approximate our system with the aid of an implicit time 
discretization and prove existence of solutions for the latter system with the help of the Leray-Schauder principle. 
Afterwards we pass to the limit in Section~\ref{sec:proof} and prove our main result on existence for the 
Cahn-Hilliard/Navier-Stokes system. Here in particular the compactness of the velocity and the attainment of the initial
data for the velocity are non-standard. Finally, in the appendix we discuss 
the corresponding results on existence of weak solutions for a similar model, which was derived in Abels, Garcke and Gr\"un 
\cite{AGG10} before.

%%%%%%%%%%%%%%%%%%%%%%%%%%%%%%%%%%%%%%%%%%%%%%%%%%%%%%%%%%%%%%%%%%%%%%%%%%%%%%%%%%%%%%%%%%%%%%%%
%%%%%%%%%%%%%%%%%%%%%%%%%%%%%%%%%%%%%%%%%%%%%%%%%%%%%%%%%%%%%%%%%%%%%%%%%%%%%%%%%%%%%%%%%%%%%%%%
%%%%%%%%%%%%%%%%%%%%%%%%%%%%%%%%%%%%%%%%%%%%%%%%%%%%%%%%%%%%%%%%%%%%%%%%%%%%%%%%%%%%%%%%%%%%%%%%
\section{Preliminaries} \label{prelimi}

%For a set  $M$ the power set will be denoted by  $\mathcal{P}(M)$ and $\chi_M$ denotes its characteristic function. 
%Moreover, we 
We denote %$\R^\dm_+=\{x\in \R^\dm: x_\dm >0\}$, 
$a\otimes b = (a_i b_j)_{i,j=1}^d$ for $a,b\in \R^d$ and $A_{\operatorname{sym}}= \frac12 (A+A^T)$ for a matrix $A\in \R^{d\times d}$.
If $X$ is a Banach space and $X'$ is its dual, then 
\begin{equation*}
  \weight{f,g} \equiv \weight{f,g}_{X',X} = f(g), \qquad f\in X', g\in X,
\end{equation*}
denotes the duality product. We write $X\hookrightarrow \hookrightarrow Y$ if $X$ is compactly embedded into $Y$. 
Moreover, if $H$ is a Hilbert space, $(\cdot\,,\cdot )_H$ denotes its inner product. 
Moreover, we use the abbreviation $(.\,,.)_{M}=(.\,,.)_{L^2(M)}$.   

\medskip

\noindent {\bf Function spaces:}
If $M\subseteq \R^d$ is measurable, 
$L^q(M)$, $1\leq q \leq \infty$ denotes the usual Lebesgue-space and $\|.\|_q$
its norm.
Moreover, $L^q(M;X)$ denotes the set of all strongly measurable
$q$-integrable functions/essentially bounded functions, where $X$ is a Banach
space. If $M=(a,b)$, we write for simplicity $L^q(a,b;X)$ and $L^q(a,b)$.
 Furthermore, $f\in L^q_{\loc}([0,\infty);X)$ if and only if $f\in L^q(0,T;X)$ for every $T>0$. 
Moreover, $L^q_{\uloc}([0,\infty); X)$ denotes the \emph{uniformly
  local} variant of $L^q(0,\infty;X)$ consisting of all strongly measurable $f\colon
[0,\infty)\to X$ such that
\begin{equation*}
  \|f\|_{L^q_{\uloc}([0,\infty); X)}= \sup_{t\geq 0}\|f\|_{L^q(t,t+1;X)} <\infty.
\end{equation*}
%If $\Omega\subseteq \R^d$ is a domain, then $f\in L^q_{\loc}(\ol{\Omega})$ if and only if 
%$f\in L^q(\Omega\cap B)$ for every ball $B$ with $B\cap \Omega \neq \emptyset$.
If $T<\infty$, we set $L^q_{\uloc}([0,T); X) := L^q(0,T;X)$.

Recall that, if $X$ is a Banach space with the Radon-Nikodym property, then 
\begin{equation*}%\label{eq:DualLq}
L^q(M;X)'= L^{q'}(M;X')\qquad \text{for every}\ 1\leq q < \infty 
\end{equation*}
by means of the duality product
$
  \weight{f,g}= \int_{M} \weight{f(x),g(x)}_{X',X}  dx 
$
for $f\in L^{q'}(M;X')$, $g\in L^{q}(M;X)$. If $X$ is reflexive or $X'$ is
separable, then $X$ has the Radon-Nikodym property, cf. Diestel and Uhl~\cite{DU77}. 

Moreover, we recall the Lemma of Aubin-Lions: If $X_0\hookrightarrow\hookrightarrow X_1 \hookrightarrow X_2$ 
are Banach spaces, $1<p<\infty$, $1\leq q <\infty$, and $I\subset \R$ is a bounded interval, then
\begin{equation}\label{eq:AubinLions}
  \left\{v\in L^p( I; X_0): \frac{dv}{dt} \in L^q(I;X_2) \right\}
  \hookrightarrow\hookrightarrow L^p(I;X_1).
\end{equation}
See J.-L.~Lions~\cite{Lio69} for the case $q>1$ and Simon~\cite{Sim87} or Roub{\'{\i}}{\v{c}}ek~\cite{Rou90} for $q=1$.

Let $\Omega \subset \R^d$ be a domain. 
Then $W^m_q(\Omega)$, $m\in \N_0$, $1\leq q\leq \infty$, denotes the usual $L^q$-Sobolev space, % $W^m_{q,\loc}(\ol{\Omega})$ its local version, 
$W^m_{q,0}(\Omega)$  the closure of $C^\infty_0(\Omega)$ in $W^m_q(\Omega)$,
$W^{-m}_q(\Omega)= (W^m_{q',0}(\Omega))'$, and $W^{-m}_{q,0}(\Omega)= (W^m_{q'}(\Omega))'$. 
% and $f\in W^{-m}_{q,\loc}(\ol{\Omega})$ if $f\in W^{-m}_q(\Omega\cap B)$ for
% every ball $B\subset \R^d$. 
The $L^2$-Bessel potential spaces are denoted by
$H^s(\Omega)$, $s\in\R$, which are defined by restriction of distributions in
$H^s(\R^d)$ to $\Omega$, cf. Triebel~\cite[Section~4.2.1]{Tri78}. We note
that, if $\Omega\subset \R^d$ is a bounded domain with $C^{0,1}$-boundary,
then there is an extension operator $E_\Omega$ which is a bounded linear operator
$E_\Omega\colon W^m_p(\Omega)\to W^m_p(\R^d)$, $1\leq p\leq \infty$ for all
$m\in \N$ and $E_\Omega f|_{\Omega} = f$ for all $f\in W^m_p(\Omega)$,
cf. Stein~\cite[Chapter VI, Section 3.2]{St70}. This extension
operator extends to $E_\Omega\colon H^s(\Omega)\to H^s(\R^d)$, which shows that
$H^s(\Omega)$ is a retract of $H^s(\R^d)$. Therefore all results on
interpolation spaces of $H^s(\R^d)$ carry over to $H^s(\Omega)$.

Given $f\in L^1(\Omega)$, we denote by 
$
  f_\Omega = \frac1{|\Omega|}\int_\Omega f(x) \,dx
$
its mean value. Moreover, for $m\in\R$ we set
\begin{equation*}
  L^q_{(m)}(\Omega):=\{f\in L^q(\Omega):f_\Omega=m\}, \qquad 1\leq q\leq \infty.  
\end{equation*}
Then for $f \in L^2(\Omega)$ we observe that 
\begin{align*}
 P_0 f:= f-f_\Omega= f-\frac1{|\Omega|}\int_\Omega f(x) \,dx
\end{align*}
is the orthogonal projection onto $L^2_{(0)}(\Omega)$. 
Furthermore,
we define
\begin{equation*}
 \Hone\equiv\Hone (\Omega)= H^1(\Omega)\cap L^2_{(0)}(\Omega), \qquad (c,d)_{\Hone(\Omega)} := (\nabla c,\nabla d)_{L^2(\Omega)}.  
\end{equation*}
Then $\Hone(\Omega)$ is a Hilbert space due to Poincar\'e's inequality.

\medskip

\noindent
{\bf Spaces of solenoidal vector-fields:}
For a bounded domain $\Omega \subset \R^d$ we denote by $C^\infty_{0,\sigma}(\Omega)$ 
in the following  the space of all
divergence free vector fields in $C^\infty_0(\Omega)^d$ and
$L^2_\sigma(\Omega)$ is its closure in the $L^2$-norm. The corresponding
Helmholtz projection is denoted by $P_\sigma$,
cf. e.g. Sohr \cite{Soh01}. We note that $P_\sigma f = f- \nabla p$, 
where $p \in W^1_2(\Omega)\cap L^2_{(0)}(\Omega)$ is the solution of the weak Neumann problem 
\begin{equation}\label{eq:WeakHelmholtz}
  (\nabla p,\nabla \varphi)_{\Omega} = (f, \nabla \varphi)\quad \text{for all}\ \varphi \in C^\infty(\ol{\Omega}).
\end{equation}

\medskip

\noindent
{\bf Spaces of continuous vector-fields:}
In the following let $I=[0,T]$ with $0<T< \infty$ or let $I=[0,\infty)$ if $T=\infty$ and let $X$ be a Banach
space. Then $BC(I;X)$ is the Banach space of all bounded and continuous
$f\colon I\to X$ equipped with the supremum norm and $BUC(I;X)$ is the
subspace of all bounded and uniformly continuous functions. Moreover, we
define $BC_w(I;X)$ as the topological vector space of all bounded and weakly
continuous functions $f\colon I\to X$. By $C^\infty_0(0,T;X)$ we denote
the vector space of all smooth functions $f\colon (0,T)\to X$ with $\supp
f\subset\subset (0,T)$.
Finally, $f\in W^1_p(0,T;X)$, $1\leq p <\infty$, if and only if $f,
\frac{df}{dt}\in L^p(0,T;X)$, where $\frac{df}{dt}$ denotes the vector-valued
distributional derivative of $f$.
Furthermore, $W^1_{p,\uloc}([0,\infty);X)$ is defined by replacing $L^p(0,T;X)$ by $L^p_{\uloc}([0,\infty);X)$ 
and we set $H^1(0,T;X)= W^1_2(0,T;X)$ and $H^1_{\uloc}([0,\infty);X) := W^1_{2,\uloc}([0,\infty);X)$.
Finally, we note:
\begin{lemma}\label{lem:CwEmbedding}
  Let $X,Y$ be two Banach spaces such that $Y\hookrightarrow X$ and $X'\hookrightarrow Y'$ densely.
  Then $L^\infty(I;Y)\cap BUC(I;X) \hookrightarrow BC_w(I;Y)$.
\end{lemma}
\noindent
For a proof, see e.g. Abels \cite{Abe09a}. 

\medskip 

\noindent
{\bf Embedding results for interpolation couples:}
Let $(X_0,X_1)$ be a compatible couple of Banach spaces, i.e., there is a Hausdorff topological 
vector space $Z$ such that $X_0,X_1 \hookrightarrow Z$, cf. Bergh and L\"ofstr\"om \cite{BL76}, 
and let $(.\,,.)_{[\theta]}$ and $(.\,,.)_{\theta,r}$, $\theta \in[0,1], 1\leq r\leq\infty$, 
denote the complex and real interpolation functor, respectively.
If in addition $X_1\hookrightarrow X_0$ densely, then for all $1\leq p<\infty$
\begin{equation}
  \label{eq:BUCEmbedding}
   W^1_p(0,T;X_0) \cap L^p(0,T;X_1) \hookrightarrow BUC(I;(X_0,X_1)_{1-\frac1p,p})
\end{equation}
continuously, cf. Amann \cite[Chapter III, Theorem 4.10.2]{Ama95}.
From this one immediately gets 
\begin{equation}
  \label{eq:BUCulocEmbedding}
   W^1_{p,\uloc}([0,T);X_0) \cap L^p_{\uloc}([0,T);X_1) \hookrightarrow BUC(I;(X_0,X_1)_{1-\frac1p,p}) \,.
\end{equation}

\medskip 

\noindent
{\bf A result related to energy inequalities:}
The following lemma will be useful for passing to the limit in energy inequalities.
\begin{lemma}\label{lem:EnergyEstim}
  Let $E\colon [0,T)\to [0,\infty)$, $0<T\leq\infty$, be a lower semi-continuous function and let $D\colon (0,T)\to [0,\infty)$ be an integrable function. Then 
    \begin{equation}\label{eq:DissEstim1}
      E(0) \varphi(0) + \int_0^T E(t) \varphi'(t) \,dt \geq \int_0^T D(t) \varphi (t)\, dt
    \end{equation}
holds for all $\varphi \in W^1_1(0,T)$ with $\varphi(T)=0$ and $\varphi \geq 0$ if and only if 
    \begin{equation}\label{eq:DissEstim2}
      E(t) + \int_s^t D(\tau) \, d\tau \leq E(s)
    \end{equation}
holds for all $s\leq t< T$ and almost all $0\leq s<T$ including $s=0$.
\end{lemma}
\noindent
For a proof, we refer to Abels \cite{Abe09a}.

%%%%%%%%%%%%%%%%%%%%%%%%%%%%%%%%%%%%%%%%%%%%%%%%%%%%%%%%%%%%%%%%%%%%%%%%%%%%%%%%%%%%%%%%%%%%%%%%%%%%%%
%%%%%%%%%%%%%%%%%%%%%%%%%%%%%%%%%%%%%%%%%%%%%%%%%%%%%%%%%%%%%%%%%%%%%%%%%%%%%%%%%%%%%%%%%%%%%%%%%%%%%%
%%%%%%%%%%%%%%%%%%%%%%%%%%%%%%%%%%%%%%%%%%%%%%%%%%%%%%%%%%%%%%%%%%%%%%%%%%%%%%%%%%%%%%%%%%%%%%%%%%%%%%
\section{Existence of weak solutions} \label{secexistence}

In this section we prove an existence result for weak solutions of the following
Navier-Stokes/ Cahn-Hilliard system for a situation with different densities. The complete
system is given by 
\begin{align}
 \partial_t (\rho \mathbf{v}) + \operatorname{div} (\rho \bv \otimes \bv) &- \operatorname{div} (2 \eta(\varphi) D \bv)
  + \nabla p + \mbox{div}(\bv \otimes \tbj) \nonumber \\
  & = - \operatorname{div}(a(\varphi) \nabla \varphi \otimes \nabla \varphi) & \mbox{in } \, Q , \label{newproblem1} \\
 \operatorname{div} \, \bv &= 0 & \mbox{in } \, Q , \label{newproblem2} \\
 \partial_t \varphi + \bv \cdot \nabla \varphi &= \mbox{div}\left(m(\varphi) \nabla \mu \right) & \mbox{in } \, Q , \label{newproblem3} \\
 \mu &= \Psi'(\varphi) + a'(\varphi) \frac{|\nabla \varphi|^2}{2} - \operatorname{div}\left( a(\varphi) \nabla \varphi \right) & \mbox{in } \, Q , \label{newproblem4} \\
 \bv|_{\partial \Omega} &= 0 & \mbox{on } \, S , \label{newproblem5} \\
 \partial_n \varphi|_{\partial \Omega} = \partial_n \mu|_{\partial \Omega} &= 0 & \mbox{on } \, S , \label{newproblem6} \\
 \left(\bv , \varphi \right)|_{t=0} &= \left( \bv_0 , \varphi_0 \right) & \mbox{in } \, \Omega , \label{newproblem7}
\end{align}
where $\tbj = -\frac{\tilde{\rho}_2 - \tilde{\rho}_1}{2} m(\varphi) \nabla \mu$. %To simplify the 
%presentation we set $\hat{\sigma} = \varepsilon = 1$, but the result will also be true for 
%general $\hat{\sigma},$ $\varepsilon > 0$.
In the above formulation and in the following, we use the abbreviations for space-time cylinders
$Q_{(s,t)}=\Omega \times (s,t)$, $Q_t = Q_{(0,t)}$ and $Q = Q_{(0,\infty)}$ and analogously for the 
boundary $S_{(s,t)}=\partial \Omega \times (s,t)$, $S_t = S_{(0,t)}$ and $S = S_{(0,\infty)}$.
Equation \eqref{newproblem5} is the no-slip boundary condition for viscous fluids, $n$ is the exterior unit normal on 
$\partial \Omega$, $\partial_n \mu |_{\partial \Omega} = 0$ means that there is no mass flux of the components 
through the boundary, and $\partial_n \p |_{\partial \Omega} = 0$ describes a contact angle of $\pi/2$ of the 
diffused interface and the boundary of the domain. 

We reformulate the first line suitably. To this end, we first calculate
\begin{align*}
 -\operatorname{div}(a(\p) \nabla \p \otimes \nabla \p) 
   = - \operatorname{div} (a(\p) \nabla \p) \nabla \p - a(\p) \nabla \left( \frac{|\nabla \p|^2}{2} \right) \,.
\end{align*}
Multiplying \eqref{newproblem4} with $\nabla \varphi$ using
$\Psi'(\p) \nabla \p = \nabla \left( \Psi(\p) \right)$ and $a'(\p) \nabla \p = \nabla \left(a(\p)\right)$
then leads to the following identity
\begin{align*}
 - \operatorname{div} (a(\p) \nabla \p) \nabla \p 
  = \mu \nabla \p - \nabla (\Psi(\p)) - \nabla \left(a(\p)\right) \frac{|\nabla \p|^2}{2} \,,
\end{align*}
which gives
\begin{align*}
 -\operatorname{div}(a(\p) \nabla \p \otimes \nabla \p) 
    = \mu \nabla \p - \nabla (\Psi(\p)) - \nabla \left( a(\p) \, \frac{|\nabla \p|^2}{2} \right) \,.
\end{align*}
With a new pressure $g = p + \Psi(\p) + a(\p) \frac{|\nabla \p|^2}{2}$ we replace \eqref{newproblem1} by 
\renewcommand{\theequation}{\arabic{section}.\arabic{equation}'}
\setcounter{equation}{0}
\begin{align} \label{replacenewproblem1}
  \partial_t(\rho \bv) + \mbox{div}(\rho \bv \otimes \bv) - \operatorname{div} (2 \eta(\p) D\bv) + \nabla g 
   + \operatorname{div}(\bv \otimes \tbj) = \mu \nabla \p \,. 
\end{align}
\renewcommand{\theequation}{\arabic{section}.\arabic{equation}}
\setcounter{equation}{7}

\subsection{Assumptions and definition of weak solutions}

In the following we summarize the assumptions needed to formulate the notion of 
a weak solution of \eqref{newproblem1}-\eqref{newproblem7} and an existence result.

\begin{assumption} \label{assumptions}
We assume that $\Omega \subset \R^d$, $d=2,3$, is a bounded domain with smooth boundary and 
additionally we impose the following conditions.
\begin{enumerate}
 \item The constitutive relation between density and phase field is given by 
       $\rho(\p) = \frac{1}{2}(\tilde{\rho}_1 + \tilde{\rho}_2) + \frac{1}{2} (\tilde{\rho}_2 - \tilde{\rho}_1) \p$
       as derived in Abels, Garcke and Gr\"un \cite{AGG11}, where $\tilde{\rho}_i>0$ are the specific constant mass 
       densities of the unmixed fluids and $\varphi$ is the difference of the volume fractions of the fluids. 
 \item We assume $a,m \in C^1(\R)$, $\eta \in C^0(\R)$ and $0 < m_0 \leq a(s),m(s),\eta(s) \leq K$ for 
       given constants $m_0,K > 0$.
 \item For the homogeneous free energy density $\Psi$ we assume that 
       $\Psi \in C([-1,1]) \cap C^2((-1,1))$ such that
       \begin{align} \label{assumptionsphi}
         \lim_{s \to -1} \Psi'(s) = -\infty \,, \quad \lim_{s \to 1} \Psi'(s) = \infty \,, \quad
         \Psi''(s) \geq -\kappa \; \mbox{ for some $\kappa \in \R$} \,.
       \end{align}
 \item Additionally we impose the condition that $\lim_{s \rightarrow \pm 1} 
       \frac{\Psi''(s)}{|\Psi'(s)|} = +\infty$. % Note: Reference on this item without labels in the next remark!
\end{enumerate}
\end{assumption}

\begin{remark} \label{remarkonassumpt}
 \begin{enumerate}
  \item An example for $\Psi$ is given by
        \begin{align} \label{examplepsi}
         \Psi(s) \ZuWeis \frac{a}{2} \left( (1+s) \,\ln(1+s) + (1-s) \, \ln(1-s) \right) - \frac{b}{2} s^2 \,, \quad s \in [-1,1] ,
        \end{align}
        where $a,b > 0$. We note that $\Psi$ is convex if and only if $a \geq b$.
  \item The Assumption \ref{assumptions} $(iv)$ is needed to reformulate the model to apply results from 
        the theory of subdifferentials.
  \item As the solution $\varphi$ will lie in the interval $[-1,1]$, we only need the functions $a,m,\eta$
        on this interval. We then extend the functions $a,m,\eta$ to the whole of $\mathbb{R}$ such that 
        $(ii)$ in Assumption \ref{assumptions} is fulfilled. 
 \end{enumerate}
\end{remark}

Now we can define a weak solution of problem \eqref{newproblem1}-\eqref{newproblem7}.
\begin{definition} \label{defweaksolution}
 Let $T \in (0,\infty]$ and set either $I = [0,\infty)$ if $T = \infty$ or $I=[0,T]$ if $T < \infty$, 
 $\bv_0 \in L^2_\sigma(\Omega)$ and $\p_0 \in H^1(\Omega)$ with $|\varphi_0| \leq 1$ almost everywhere in $\Omega$.
 If in addition Assumption \ref{assumptions} holds we call the triple $(\bv,\p,\mu)$ 
 with the properties
 \begin{align*}
  & \bv \in BC_w(I;L^2_\sigma(\Omega)) \cap L^2(0,T;H_0^1(\Omega)^d) \,, \\
  & \p \in BC_w(I;H^1(\Omega)) \cap L^2_{\mbox{\footnotesize uloc}}(I;H^2(\Omega)) \,, \; \
        \Psi'(\p) \in L^2_{\mbox{\footnotesize uloc}}(I;L^2(\Omega)) \,, \\
  & \mu \in L^2_{\mbox{\footnotesize uloc}}(I;H^1(\Omega)) 
    \; \mbox{ with } \; \nabla \mu \in L^2(0,T;L^2(\Omega)) 
 \end{align*}
 a weak solution of \eqref{newproblem1}-\eqref{newproblem7},
 if the following conditions are satisfied.
 \begin{align}  
   - \left(\rho \bv , \partial_t \bpsi \right)_{Q_T} 
  &+ \left( \operatorname{div}(\rho \bv \otimes \bv) , \bpsi \right)_{Q_T}
  + \left(2 \eta(\p) D\bv , D\bpsi \right)_{Q_T} 
   - \left( (\bv \otimes \tbj) , \nabla \bpsi \right)_{Q_T} \nonumber \\
  &= \left( \mu \nabla \varphi , \bpsi \right)_{Q_T} \label{weakline1} 
 \end{align}
 for all $\bpsi \in \left[C_0^\infty(\Omega \times (0,T))\right]^d$ with $\operatorname{div} \bpsi = 0$,
 \begin{align} 
  - \left(\p , \partial_t \zeta \right)_{Q_T} 
  + \left( \bv \cdot \nabla \p , \zeta \right)_{Q_T}
  &= - \left(m(\p) \nabla \mu , \nabla \zeta \right)_{Q_T}  \label{weakline2} 
 \end{align}
 for all $\zeta \in C_0^\infty((0,T);C^1(\overline{\Omega}))$,
 \begin{align}
   \mu = \Psi'(\varphi) + a'(\varphi) \frac{|\nabla \varphi|^2}{2} &- \operatorname{div}\left(a(\varphi) \nabla \varphi \right)
     \; \mbox{ almost everywhere in } \, Q_T \; \mbox{and}  \label{weakline3} \\
  \left.\left( \bv,\p \right)\right|_{t=0} &= \left( \bv_0 , \p_0 \right) \,. \label{weakline4} 
 \end{align}
 Moreover,
 \begin{align} 
  E_{\mbox{\footnotesize tot}}(\p(t),\bv(t)) &+ \int_{Q_{(s,t)}} 2 \eta(\p) \, |D\bv|^2 \, d(x,\tau) 
        + \int_{Q_{(s,t)}} m(\p) |\nabla \mu|^2 \, d(x,\tau) \nonumber \\
   &\leq E_{\mbox{\footnotesize tot}}(\p(s),\bv(s)) \label{weakline5}
 \end{align}
 for all $t \in [s,\infty)$ and almost all $s \in [0,\infty)$ has to hold (including $s=0$).
 The total energy $E_{\mbox{\footnotesize tot}}$ is the sum of the kinetic and the free energy, see \eqref{totalenergy}.
\end{definition}

%[Remark: We can also take $\left(v \cdot \nabla \p,\eta \right)_{Q_T}$ instead of 
%$-\left( \varphi v , \nabla \eta \right)_{Q_T}$ in \eqref{weakline2}.]
\subsection{Existence theorem}

Our main result of this work is the following existence theorem for weak solutions. 

\begin{theorem} \label{existenceweaksolution}
 Let Assumption \ref{assumptions} hold, $\bv_0 \in L^2_\sigma(\Omega)$, $\varphi_0 \in H^1(\Omega)$
 with $|\varphi_0| \leq 1$ almost everywhere and $\int_\Omega \hspace*{-14pt}{-}\hspace*{5pt} \varphi_0 \, dx \in (-1,1)$. 
 Then there exists a weak solution $(\bv,\varphi,\mu)$ of \eqref{newproblem1}-\eqref{newproblem7} 
 in the sense of Definition \ref{defweaksolution}.
\end{theorem}

In order to prove the theorem, we reformulate line \eqref{weakline3} to an equivalent equation.
Therefore we introduce the function $A(s) \ZuWeis \int_0^s \sqrt{a(\tau)} \, d\tau$.
Then $A'(s) = \sqrt{a(s)}$ and
\begin{align*}
 - \sqrt{a(\p)} \, \Delta A(\p) = a'(\p) \, \frac{|\nabla \p|^2}{2} - \mbox{div} \left(a(\p) \, \nabla \p \right)
\end{align*}
resulting from a straightforward calculation.
\iffalse
\begin{align*}
 - \sqrt{a(\p)} \, \Delta A(\p) =& - \sqrt{a(\p)} \, \mbox{div} \left( \nabla A(\p) \right)
  = - \sqrt{a(\p)} \, \mbox{div} \left( A'(\p) \, \nabla \p \right)
  = - \sqrt{a(\p)} \, \mbox{div} \left( \sqrt{a(\p)} \, \nabla \p \right) \\
  =& - \sqrt{a(\p)} \left( \nabla \left( \sqrt{a(\p)} \right) \cdot \nabla \p + \sqrt{a(\p)} \Delta \p \right) \\
  =& - \sqrt{a(\p)} \frac{1}{2} a(\p)^{-\frac{1}{2}} \, a'(\p) \, \nabla \p \cdot \nabla \p - a(\p) \, \Delta \p \\
  =& - a'(\p) \, \frac{|\nabla \p|^2}{2} - a(\p) \, \Delta \p \\
  =& \; a'(\p) \, \frac{|\nabla \p|^2}{2} - \underbrace{a'(\p) \, \nabla \p}_{=\nabla a(\p)} \cdot \nabla \p - a(\p) \, \Delta \p \\
  =& \; a'(\p) \, \frac{|\nabla \p|^2}{2} - \mbox{div} \left(a(\p) \, \nabla \p \right) \,. 
\end{align*}
\fi

With this notation \eqref{weakline3} reduces to
\begin{align} \label{replaceInewproblem4}
 \mu = \Psi'(\varphi) - \sqrt{a(\p)} \, \Delta A(\p) \,.
\end{align}
We also rewrite the free energy with the help of $A$ to 
\begin{align*}
 E_{\mbox{\footnotesize free}}(\p) = \int_\Omega \left( \Psi(\p) + \frac{|\nabla A(\p)|^2}{2} \right) dx \,.
\end{align*}
The kinetic energy is given by $E_{\mbox{\footnotesize kin}}(\p,v) = \int_\Omega \rho \frac{|v|^2}{2} \, dx$ and the 
total energy as the sum of the kinetic and free energy
\begin{align} \label{totalenergy}
 E_{\mbox{\footnotesize tot}}(\p,v)  
         = E_{\mbox{\footnotesize kin}}(\p,v) + E_{\mbox{\footnotesize free}}(\p)
         = \int_\Omega \rho \frac{|v|^2}{2} \, dx + \int_\Omega \left( \Psi(\p) + \frac{|\nabla A(\p)|^2}{2} \right) dx .
\end{align}
The next step is to rewrite \eqref{replaceInewproblem4} with the help of a subdifferential.
To this end, we set $[a,b] \ZuWeis A([-1,1])$ and define a reparametrized potential $\widetilde{\Psi}$ through
\begin{align} \label{deftildepsi}
 \widetilde{\Psi} : \R \to \R \,, \quad 
   \widetilde{\Psi}(r) \ZuWeis \left\{ \begin{array}{cl} 
                                    \Psi(A^{-1}(r)) & \mbox{if } \, r \in [a,b] ,\\
                                    + \infty   & \mbox{else} \,.
                                   \end{array}
                               \right.
\end{align}
 This $\widetilde{\Psi}$ fulfills analogous assumptions as $\Psi$, that is
 $\widetilde{\Psi} \in C^0([a,b]) \cap C^2((a,b))$ and 
 \begin{align} \label{assumptildepsi}
   \lim_{r \to a} \widetilde{\Psi}'(r) = - \infty \,, \quad 
   \lim_{r \to b} \widetilde{\Psi}'(r) = + \infty \,, \quad
   \widetilde{\Psi}''(r) \geq - \widetilde{\kappa} \; \mbox{ for all } \; r \in (a,b) 
 \end{align}
for some $\widetilde{\kappa} \in \mathbb{R}$. Here Assumption \ref{assumptions} $(iv)$ is needed.
We define $\widetilde{\Psi}_0(r) \ZuWeis \widetilde{\Psi}(r) + \frac{\widetilde{\kappa}}{2} r^2$
to get a convex function $\widetilde{\Psi}_0$. In particular it holds 
$\widetilde{\Psi}'(r) = \widetilde{\Psi}_0'(r) - \widetilde{\kappa} r$ and 
$\widetilde{\Psi}'(A(s)) = \widetilde{\Psi}_0'(A(s)) - \widetilde{\kappa} A(s)$ for $s \in [-1,1]$. Furthermore we observe
$\widetilde{\Psi}(A(s)) = \Psi(s)$ for $s \in [-1,1]$ and therefore 
$\Psi'(s) = \widetilde{\Psi}'(A(s)) \cdot \sqrt{a(s)}$, in particular $\left(a(s)\right)^{-\frac{1}{2}} \cdot \Psi'(s) = \widetilde{\Psi}'(A(s))$.

This notation leads to a reformulation of \eqref{replaceInewproblem4} to
\begin{align} \label{replaceIInewproblem4}
 \left(a(\varphi)\right)^{-\frac{1}{2}} \mu + \widetilde{\kappa} A(\varphi) = \widetilde{\Psi}_0'(A(\varphi)) - \Delta A(\varphi) \,.
\end{align}
%We note that it would be possible to choose w.l.o.g. $\tilde{\kappa} = \kappa$ and $\tilde{\kappa} \geq 0$.
Now we use a result from Abels and Wilke \cite{AW07} for the energy $\widetilde{E} : L^2(\Omega) \to \R$ 
with domain $\mbox{dom}\, \widetilde{E} = \{u \in H^1(\Omega) \;|\; a \leq u \leq b \,\mbox{ a.e.} \}$ given by 
\begin{align} \label{energytilde}
 \widetilde{E}(u) = \left\{ 
                    \begin{array}{cl} 
                     \frac{1}{2} \int_\Omega |\nabla u|^2 \, dx + \int_\Omega \widetilde{\Psi}_0(u) \, dx
                        & \mbox{for } \; u \in \mbox{dom}\, \widetilde{E} \,, \\
                     + \infty & \mbox{else} \,.
                    \end{array}
                    \right.
\end{align}
From Theorem 3.12.8 in \cite{Abe07}, which is a variant of Theorem 4.3 in \cite{AW07}, the domain of definition of the 
subgradient $\partial \widetilde{E}$ is given by
\begin{align*}
 \D(\partial \widetilde{E}) = \{ u \in H^2(\Omega) \;|\; \widetilde{\Psi}_0'(u) \in L^2(\Omega) \,, \;
            \widetilde{\Psi}_0''(u) |\nabla u|^2 \in L^1(\Omega) \,, \; \partial_n u |_{\partial \Omega} = 0 \}
\end{align*}
and for $u \in \D(\partial \widetilde{E})$ it holds that $\partial \widetilde{E}(u) = - \Delta u + \widetilde{\Psi}_0'(u)$.
Furthermore there holds the estimate
\begin{align} \label{AWestimate}
 \|u\|^2_{H^2} + \|\widetilde{\Psi}_0'(u)\|^2_{L^2} + \int_\Omega \widetilde{\Psi}_0''(u(x)) |\nabla u(x)|^2 \, dx \leq 
  C \left( \| \partial \widetilde{E}(u) \|_{L^2}^2 + \|u\|_{L^2}^2 + 1 \right) \,.
\end{align}
Analogously we define the energy $E : L^2(\Omega) \to \R$ 
with domain $\mbox{dom}\, E = \{\p \in H^1(\Omega) \;|\; -1 \leq \p \leq 1 \,\mbox{ a.e.} \}$ given by 
\begin{align} \label{helpenergy}
 E(\p) = \left\{ 
                    \begin{array}{cl} 
                     \frac{1}{2} \int_\Omega |\nabla \p|^2 \, dx + \int_\Omega \Psi_0(\p) \, dx
                        & \mbox{for } \; \varphi \in \mbox{dom}\, E \,, \\
                     + \infty & \mbox{else} \,.
                    \end{array}
                    \right.
\end{align}
Here it holds that $\D(\partial E) = \{ \p \in H^2(\Omega) \;|\; \Psi_0'(\p) \in L^2(\Omega) \,, \;
            \Psi_0''(\p) |\nabla \p|^2 \in L^1(\Omega) \,, \; \partial_n \p |_{\partial \Omega} = 0 \} $.
One can show that that $A(\p) \in \mathcal{D}(\partial \widetilde{E})$ if and only if $\p \in \mathcal{D}(\partial E)$
and therefore we get for $u=A(\p)$ with $\p \in \D(\partial E)$ the identity
\begin{align} \label{subgradientenergytilde}
 \partial \widetilde{E}(A(\p)) = - \Delta A(\p) + \widetilde{\Psi}_0'(A(\p)) \,.
\end{align}
This leads finally to the reformulation of \eqref{weakline3} as
\begin{align} \label{replaceIIInewproblem4}
 \left(a(\varphi)\right)^{-\frac{1}{2}} \mu + \widetilde{\kappa} A(\varphi) = \partial \widetilde{E}(A(\varphi))  \,.
\end{align}
The fact that the right hand side equals the subgradient of $\widetilde{E}$ at the point $A(\p)$
will be of importance in the following analysis.

With the above notation we can rewrite the free energy as 
\begin{align*}
 E_{\mbox{\footnotesize free}}(\p) &= \int_\Omega \Psi(\p) \, dx + \int_\Omega \frac{|\nabla A(\p)|^2}{2} \, dx
  = \int_\Omega \widetilde{\Psi}(A(\p)) \, dx + \int_\Omega \frac{|\nabla A(\p)|^2}{2} \, dx \\
  &= \int_\Omega \widetilde{\Psi}_0(A(\p)) \, dx + \int_\Omega \frac{|\nabla A(\p)|^2}{2} \, dx - \int_\Omega \frac{\tilde{\kappa}}{2} (A(\p))^2 \, dx \\
  &= \widetilde{E}(A(\p)) - \frac{\tilde{\kappa}}{2} \|A(\p)\|_{L^2}^2 \,.
\end{align*}
We summarize the reformulation of problem \eqref{newproblem1}-\eqref{newproblem7}:
\begin{align}
 \partial_t (\rho \mathbf{v}) + \mbox{div} (\rho \bv \otimes \bv) &- \mbox{div} (2 \eta(\varphi) D \bv)
  + \nabla g + \mbox{div}(\bv \otimes \tbj) = \mu \nabla \p & \mbox{in } \, Q , \label{equivnewproblem1} \\
 \mbox{div} \, \bv &= 0 & \mbox{in } \, Q , \label{equivnewproblem2} \\
 \partial_t \varphi + \bv \cdot \nabla \varphi &= \mbox{div}\left(m(\varphi) \nabla \mu \right) & \mbox{in } \, Q , \label{equivnewproblem3} \\
 \left(a(\varphi)\right)^{-\frac{1}{2}} \mu + \widetilde{\kappa} A(\varphi) &= \widetilde{\Psi}_0'(A(\varphi)) - \Delta A(\varphi) & \mbox{in } \, Q , \label{equivnewproblem4} \\
 \bv|_{\partial \Omega} &= 0 & \mbox{on } \, S , \label{equivnewproblem5} \\
 \partial_n \varphi|_{\partial \Omega} = \partial_n \mu|_{\partial \Omega} &= 0 & \mbox{on } \, S , \label{equivnewproblem6} \\
 \left(\bv , \varphi \right)|_{t=0} &= \left( \bv_0 , \varphi_0 \right) & \mbox{in } \, \Omega , \label{equivnewproblem7}
\end{align}
where $\tbj = -\frac{\tilde{\rho}_2 - \tilde{\rho}_1}{2} m(\varphi) \nabla \mu$.

\section{Implicit time discretization}\label{sec:Implicit}

\subsection{Definition of the time-discrete problem} \label{subsec:DefImplicit}

In order to prove Theorem \ref{existenceweaksolution}, we use an implicit time discretization. 
To this end, let $h=\frac{1}{N}$ for $N \in \mathbb{N}$ and $\bv_k \in L^2_\sigma(\Omega)$, 
$\p_k \in H^1(\Omega)$ with $\Psi'(\p_k) \in L^2(\Omega)$ and $\rho_k = 
\frac{1}{2}(\tilde{\rho}_1 + \tilde{\rho}_2) + \frac{1}{2} (\tilde{\rho}_2 - \tilde{\rho}_1) \p_k$ 
be given. We construct $(\bv,\p,\mu)=(\bv_{k+1},\p_{k+1},\mu_{k+1})$ as solution of the 
following non-linear system, where
\begin{align*}
 \tbj = \tbj_{k+1} = - \tfrac{\tilde{\rho}_2 - \tilde{\rho}_1}{2} m(\varphi_k) \nabla \mu_{k+1}
  = - \tfrac{\tilde{\rho}_2 - \tilde{\rho}_1}{2} m(\varphi_k) \nabla \mu \,.
\end{align*}
Find $(\bv,\p,\mu)$ with $\bv \in H_0^1(\Omega)^d \cap L^2_\sigma(\Omega)$, $\varphi \in \mathcal{D}(\partial E)$
and $\mu \in H^2_n(\Omega) = \{u \in H^2(\Omega) \,|\, \left. \partial_n u \right|_{\partial \Omega} = 0 \mbox{ on } \partial \Omega\}$, such that
\begin{align} 
 \left( \frac{\rho \bv - \rho_k \bv_k}{h} , \bpsi\right)_\Omega 
 &+ \left( \mbox{div}(\rho_k  \bv \otimes \bv) , \bpsi \right)_{\Omega}
 + \left(2 \eta(\p_k) D\bv , D \bpsi \right)_\Omega 
 + \left( \mbox{div}( \bv \otimes \tbj) , \bpsi \right)_{\Omega} \nonumber \\
 &= \left( \mu \nabla \p_k , \bpsi \right)_\Omega \label{timediscretizationline1}
\end{align}
for all $\bpsi \in C_{0,\sigma}^\infty(\Omega)$, 
%(or equivalently $\bpsi \in H^1_0(\Omega)^d \cap L^2_\sigma(\Omega)$), 
\begin{align} 
 &\frac{\p - \p_k}{h} + \bv \cdot \nabla \p_k = \mbox{div} \left( m(\p_k) \nabla \mu \right) 
  \; \mbox{ almost everywhere in $\Omega$ and}  \label{timediscretizationline2} \\
 & \frac{\p - \p_k}{A(\p) - A(\p_k)} \mu + \widetilde{\kappa} \, \frac{A(\p) + A(\p_k)}{2} 
 = -\Delta A(\p) + \widetilde{\Psi}_0'(A(\p)) \; \mbox{ almost everywhere in $\Omega$.} \label{timediscretizationline3}
\end{align}

\vskip 10pt
\begin{remark} %\vskip -10pt
 \begin{enumerate}
  \item Multiplying identity \eqref{timediscretizationline2} with $-\frac{\tilde{\rho}_2 - \tilde{\rho}_1}{2}$ leads to
        \begin{align*}
         -\frac{\rho - \rho_k}{h} - \bv \cdot \nabla \rho_k = \operatorname{div} \tbj \,,
        \end{align*}
        where we used the definition of $\tbj = \tbj_{k+1}$ and the linear dependence
        $\rho(\p) = \frac{1}{2}(\tilde{\rho}_1 + \tilde{\rho}_2) + \frac{1}{2} (\tilde{\rho}_2 - \tilde{\rho}_1) \p$.
        Using $\operatorname{div}(\bv \otimes \tbj) = (\operatorname{div} \tbj) \bv 
        + \left(\tbj \cdot \nabla \right) \bv$ leads to an equivalent version of \eqref{timediscretizationline1}
        given by
        \begin{align} 
         &\left( \frac{\rho \bv - \rho_k \bv_k}{h} , \bpsi\right)_\Omega 
          + \left( \operatorname{div}(\rho_k  \bv \otimes \bv) , \bpsi \right)_{\Omega}
          + \left(2 \eta(\p_k) D\bv , D \bpsi \right)_\Omega \nonumber \\
          +& \left( \left(\operatorname{div} \tbj - \frac{\rho - \rho_k}{h} - \bv \cdot \nabla \rho_k \right) \frac{\bv}{2} 
              , \bpsi \right)_{\Omega} 
          + \left( \left( \tbj \cdot \nabla \right) \bv , \bpsi \right)_{\Omega} 
          = \left( \mu \nabla \p_k , \bpsi \right)_\Omega  \label{equivtimediscretizationline1}
        \end{align}
        for all $\bpsi \in C_{0,\sigma}^\infty(\Omega)$. We will use this equivalent version in the 
        following especially in the a-priori estimate for solutions of the time-discrete problem. 
  \item Integrating equation \eqref{timediscretizationline2} with respect to the spatial variable, using $\operatorname{div} \, \bv = 0$ and the boundary
        conditions, we obtain $\int_\Omega \varphi \, dx = \int_\Omega \varphi_k \, dx$, which means that 
        $\int_\Omega \varphi_k \, dx = \int_\Omega \varphi_0 \, dx$ is constant.
 \end{enumerate}
\end{remark}

\begin{lemma} \label{derivativeforchempot}
 Assume that $\varphi \in \mathcal{D}(\partial E)$ and 
 $\mu \in H^1(\Omega)$ solve \eqref{timediscretizationline3} with given $\varphi_k \in H^2(\Omega)$ 
 with $|\varphi_k| \leq 1$ in $\Omega$ such that
 \begin{align*}
  \tfrac{1}{|\Omega|} \int_\Omega \varphi \, dx = \tfrac{1}{|\Omega|} \int_\Omega \varphi_k \, dx \in (-1,1) \,.
 \end{align*}
 Then there is a constant $C=C(\int_\Omega \varphi_k)>0$, such that
 \begin{align*}
  \| \widetilde{\Psi}_0'(A(\p)) \|_{L^2(\Omega)} + \left| \int_\Omega \mu \, dx \right| 
   & \leq C (\|\nabla \mu\|_{L^2} + \|\nabla \varphi\|_{L^2}^2 + \|\nabla \varphi_k\|_{L^2}^2  + 1) \; \mbox{and} \\
  \|\partial \widetilde{E}(A(\p)) \|_{L^2(\Omega)} & \leq C \left( \|\mu\|_{L^2(\Omega)} + 1 \right) \,. 
 \end{align*}
 %if $\|(\p,\p_k)\|_{H^1} \leq R$.
\end{lemma}
\begin{proof}
First we note that $F(\p,\p_k) := \frac{A(\p)-A(\p_k)}{\p - \p_k}$ fulfills $0 < c \leq F(\p,\p_k) \leq C$ 
for some $c,C$ independent of $k$ due to
\begin{align*}
 \frac{A(\p)-A(\p_k)}{\p - \p_k} = \int_0^1 A'(\tau \p + (1-\tau) \p_k) \, d\tau 
   = \int_0^1 \sqrt{a(\tau \p + (1-\tau)\p_k)} \, d\tau
\end{align*}
and $0 < m_0 \leq a$, $a \in C^1(\mathbb{R})$ and $\p(x) \in[-1,1]$ almost everywhere. 
Testing \eqref{timediscretizationline3} with $\zeta = F(\p,\p_k) (\p - \overline{\p})$, 
where $\overline{\varphi} = \frac{1}{|\Omega|} \int_\Omega \varphi \, dx$ is the mean value of $\varphi$ we get
\begin{align} \label{testedwithFanddiff}
  &\int_\Omega \mu (\p-\overline{\p} ) \, dx
 + \int_\Omega \tilde{\kappa} \frac{A(\p) + A(\p_k)}{2} F(\p,\p_k) (\p - \overline{\p}) \, dx \nonumber \\
 =& \int_{\Omega} - \Delta A(\p) F(\p,\p_k) (\p - \overline{\p}) \, dx 
   + \int_\Omega \widetilde{\Psi}_0'(A(\p)) F(\p,\p_k) (\p - \overline{\p}) \, dx \,.
\end{align}
With $\mu_0 = \mu - \overline{\mu}$ we observe that $\int_\Omega \mu (\p - \overline{\p}) \, dx = \int_\Omega \mu_0 \p \, dx$.
Furthermore, due to
\begin{align} \label{helpestforquot}
 \nabla& \left( \frac{A(\p)-A(\p_k)}{\p - \p_k} \right) \nonumber \\
  &= \int_0^1 \left( a^{-\frac{1}{2}}(\tau \p + (1-\tau)\p_k) \, a'(\tau \p + (1-\tau) \p_k) \, (\tau \nabla \p + (1-\tau) \nabla \p_k) \right) d\tau \,,
\end{align}
the boundedness of $a \geq m_0 > 0$ and $a'$ on $[-1,1]$ and $\left(\tau \p + (1-\tau) \p_k\right) \in [-1,1]$, we get
\begin{align*}
 \left| \nabla \left( \frac{A(\p)-A(\p_k)}{\p - \p_k} \right) \right| \leq C \left( |\nabla \p| + |\nabla \p_k| \right) ,
\end{align*}
which leads to
\begin{align*}
 \left| \int_{\Omega} - \Delta A(\p) \cdot F(\p,\p_k) (\p - \overline{\p}) \, dx \right|
   &= \left| \int_\Omega \underbrace{\nabla A(\p)}_{= \sqrt{a(\p)} \nabla \p} 
         \cdot \nabla \left( \left( \frac{A(\p)-A(\p_k)}{\p - \p_k} \right) \, (\p - \overline{\p}) \right) dx \right| \\
  & \leq C \int_\Omega |\nabla \p|^2 \, dx + C \int_\Omega |\nabla \p| \left( |\nabla \p| + |\nabla \p_k| \right) dx \\
  & \leq C \left( \|\nabla \p\|_{L^2}^2 + \|\nabla \p_k\|_{L^2}^2 \right) \,.
\end{align*}
To estimate the last integral in \eqref{testedwithFanddiff} we calculate at first
$\widetilde{\Psi}_0'(A(\varphi)) = \Psi_0'(\varphi) + \widetilde{\kappa} A(\varphi) - \kappa \varphi$
and use the fact that $\overline{\p} \in (-1+\varepsilon,1-\varepsilon)$ for some $\varepsilon > 0$. In 
addition with the assumption $\lim_{\p \to \pm1} \Psi_0'(\p) = \pm \infty$ one can then show the inequality
$\Psi_0'(\p) (\p - \overline{\p}) \geq C_\varepsilon |\Psi_0'(\p)| - \tilde{C}_\varepsilon$ in three steps 
in the intervals $[-1,-1 + \frac{\varepsilon}{2}]$, $[-1+\frac{\varepsilon}{2},1-\frac{\varepsilon}{2}]$ and
$[1-\frac{\varepsilon}{2},1]$ successively. Altogether this leads to an estimate of the following form:
\begin{align*}
 \int_\Omega \widetilde{\Psi}_0'(A(\p)) F(\p,\p_k) (\p - \overline{\p}) \, dx
   \geq C \int_\Omega |\Psi_0'(\p)| \, dx - C_1 \,.
\end{align*}
Using the last inequalities and \eqref{testedwithFanddiff} we get
\begin{align*}
 \int_\Omega |\Psi_0'(\p)| \, dx
  \leq& \; C \|\mu_0\|_{L^2} \|\p\|_{L^2}  
   + C \int_\Omega \frac{\tilde{\kappa}}{2} 
         \underbrace{\left|A(\p)+A(\p_k)\right|}_{\leq C(|\p|+|\p_k|)} 
           \left|\frac{A(\p)-A(\p_k)}{\p-\p_k}\right| \left|\p-\overline{\p}\right| dx \\
  & \; + C(\|\nabla \p\|_{L^2}^2 + \|\nabla \p_k\|_{L^2}^2) + C_1 \\
  \leq & \; C (\|\mu_0\|_{L^2} + \|\nabla \varphi\|_{L^2}^2 + \|\nabla \varphi_k\|_{L^2}^2  + 1) \\
  \leq & \; C (\|\nabla \mu\|_{L^2} + \|\nabla \varphi\|_{L^2}^2 + \|\nabla \varphi_k\|_{L^2}^2  + 1)  \,.
  %\leq & \; C(\|\p\|_{H^1},\|\p_k\|_{H^1}) \, (\|\mu_0\|_{L^2} + 1) 
  %\; \leq \; C(\|\p\|_{H^1},\|\p_k\|_{H^1}) \, (\|\nabla \mu\|_{L^2} + 1) \,.
\end{align*}
where we have use $|A(\p)| \leq C |\p|$ and the bounds $|\varphi|$, $|\varphi_k| \leq 1$. 
Now we test \eqref{timediscretizationline3} with $\frac{A(\p)-A(\p_k)}{\p-\p_k}$ to get
\begin{align*}
 \int_\Omega \mu \, dx = &\int_\Omega \nabla A(\p) \cdot \nabla \left(\frac{A(\p)-A(\p_k)}{\p-\p_k}\right) dx
  - \int_\Omega \widetilde{\Psi}_0'(A(\p)) \cdot \frac{A(\p)-A(\p_k)}{\p-\p_k} \, dx \\
  & - \int_\Omega \frac{\tilde{\kappa}}{2} \left(A(\p) + A(\p_k)\right) \cdot \frac{A(\p)-A(\p_k)}{\p-\p_k} \, dx \,,
\end{align*}
which together with the previous estimates leads to 
\begin{align*}
 \left| \int_\Omega \mu \, dx \right| \leq & 
   C (\|\nabla \mu\|_{L^2} + \|\nabla \varphi\|_{L^2}^2 + \|\nabla \varphi_k\|_{L^2}^2  + 1) \,.
\end{align*}
Finally, the estimates of the subdifferential $\partial \widetilde{E}(A(\varphi)) 
= - \Delta A(\varphi) + \widetilde{\Psi}_0'(A(\varphi))$ and of $\widetilde{\Psi}_0'(A(\varphi)) = \Psi_0'(\varphi)$,
each in the $L^2(\Omega)$-norm,
follow directly from \eqref{timediscretizationline3} and inequality \eqref{AWestimate}.
\end{proof}

\subsection{Existence of a solution of the time-discrete problem}

\begin{lemma} \label{timediscreteexistence}
 Let $v_k \in L^2_\sigma(\Omega)$, $\varphi_k \in H^2(\Omega)$ %with $\Psi'(\varphi_k) \in L^2(\Omega)$ 
 and $\rho_k = \frac{1}{2} (\tilde{\rho}_1 + \tilde{\rho}_2) + \frac{1}{2} (\tilde{\rho}_2 - \tilde{\rho}_1) \varphi_k$
 be given. Then there are some $(\bv,\p,\mu) \in \left( H^1_0(\Omega)^d \cap L^2_\sigma(\Omega) \right)
 \times \mathcal{D}(\partial E) \times H^2_n(\Omega)$ solving \eqref{timediscretizationline1}-\eqref{timediscretizationline3},
 which satisfy in addition the discrete energy estimate
 \begin{align} \label{discreteenergyestimate}
  E_{\mbox{\footnotesize tot}}&(\p,\bv) + \int_\Omega \rho_k \frac{|\bv - \bv_k|^2}{2} \, dx
    + \int_\Omega \frac{|\nabla A(\p) - \nabla A(\p_k)|^2}{2} \, dx \nonumber \\
    &+ h \int_\Omega 2 \eta(\p_k) |D\bv|^2 \, dx + h \int_\Omega m(\p_k) |\nabla \mu|^2 \, dx 
    \leq E_{\mbox{\footnotesize tot}}(\p_k,\bv_k) \,. 
 \end{align}
\end{lemma}
\begin{proof}
First we show the a-priori estimate \eqref{discreteenergyestimate} for any 
$(\bv,\p,\mu) \in \left( H^1_0(\Omega)^d \cap L^2_\sigma(\Omega) \right) \times \mathcal{D}(\partial E) \times H^2_n(\Omega)$ 
solving \eqref{timediscretizationline1}-\eqref{timediscretizationline3}.

In order to test \eqref{timediscretizationline1} with $\bpsi = \bv$ we need some preparations.
First we observe
\begin{align*}
 \int_\Omega \left( (\mbox{div} \, \tbj) \frac{\bv}{2} + \left( \tbj \cdot \nabla \right) \bv \right) \cdot \bv \, dx
  = \int_\Omega \mbox{div} \left( \tbj \frac{|\bv|^2}{2} \right) \, dx = 0 
\end{align*}
and then we calculate
\begin{align*}
 \int_\Omega & \left( \mbox{div}(\rho_k \bv \otimes \bv) - (\nabla \rho_k \cdot \bv) \frac{\bv}{2} \right) \cdot \bv \, dx \\
 &= \int_\Omega \left( \mbox{div}(\rho_k \bv \otimes \bv) - \mbox{div}(\rho_k \bv) \frac{\bv}{2} \right) \cdot \bv \, dx \\
 &= \int_\Omega \left( \mbox{div} (\rho_k \bv) |\bv|^2 
      + \rho_k \bv \cdot \nabla \left( \frac{|\bv|^2}{2} \right) 
      - \mbox{div}(\rho_k \bv) \frac{|\bv|^2}{2} \right) dx\\
 &= \int_\Omega \mbox{div} \left( \rho_k \bv \frac{|\bv|^2}{2} \right) \, dx 
 = 0 \,,
\end{align*}
where we used integration by parts and $\operatorname{div} \bv = 0$. Furthermore from the simple algebraic equation
\begin{align*}
 \mathbf{a} \cdot (\mathbf{a}-\mathbf{b}) = \frac{|\mathbf{a}|^2}{2} - \frac{|\mathbf{b}|^2}{2} + \frac{|\mathbf{a}-\mathbf{b}|^2}{2} 
  \quad \mbox{ for } \; \mathbf{a},\mathbf{b} \in \R^d
\end{align*}
we get 
\begin{align*}
 \frac{1}{h} \left(\rho \bv - \rho_k \bv_k \right) \cdot \bv 
  &= \frac{1}{h} (\rho - \rho_k)\,|\bv|^2 + \frac{1}{h} \rho_k \, (\bv-\bv_k)\cdot \bv \\
  &= \frac{1}{h} (\rho - \rho_k)\,|\bv|^2 + \frac{1}{h} \rho_k \left( \frac{|\bv|^2}{2} - \frac{|\bv_k|^2}{2} \right) 
                 + \frac{1}{h} \rho_k \frac{|\bv-\bv_k|^2}{2} \\
  &= \frac{1}{h} \left( \rho \frac{|\bv|^2}{2} - \rho_k \frac{|\bv_k|^2}{2} \right) 
        + \frac{1}{h} \, (\rho - \rho_k) \, \frac{|\bv|^2}{2} 
        + \frac{1}{h} \rho_k \frac{|\bv-\bv_k|^2}{2} \,.
\end{align*}
Testing \eqref{timediscretizationline1} or equivalently \eqref{equivtimediscretizationline1} with $\bpsi = \bv$
and using the above identities we obtain
\begin{align} \label{testline1}
 0 = \int_\Omega \frac{\rho |\bv|^2 - \rho_k |\bv_k|^2}{2h} \, dx
 + \int_\Omega \rho_k \frac{|\bv-\bv_k|^2}{2h} \, dx
 + \int_\Omega 2 \eta(\p_k) |D\bv|^2 \, dx
 - \int_\Omega \mu \, (\nabla \p_k \cdot \bv) \, dx \,.
\end{align}
Moreover, choosing $\mu$ in \eqref{timediscretizationline2} as a test function, we get
\begin{align} \label{testline2}
 0 = \int_\Omega \frac{\p - \p_k}{h} \, \mu \, dx
  + \int_\Omega (\bv \cdot \nabla \p_k) \, \mu \, dx
  + \int_\Omega m(\p_k) |\nabla \mu|^2 \, dx \,.
\end{align}
Finally we test \eqref{timediscretizationline3} with $\frac{1}{h} (A(\p) - A(\p_k))$ to get
\begin{align} \label{testline3}
 0 = & \frac{1}{h} \int_\Omega \nabla A(\p) \cdot \nabla (A(\p) - A(\p_k)) \, dx
        + \int_\Omega \widetilde{\Psi}_0'(A(\p)) \frac{A(\p)-A(\p_k)}{h} \, dx \nonumber \\
     & - \int_\Omega \mu \, \frac{\p - \p_k}{h} \, dx
      - \int_\Omega \widetilde{\kappa} \frac{A(\p)^2 - A(\p_k)^2}{2h} \, dx \,.
\end{align}
Summing the identities \eqref{testline1}-\eqref{testline3} leads to
\begin{align*}
 0  =& \int_\Omega \frac{\rho |\bv|^2 - \rho_k |\bv_k|^2}{2h} \, dx
     + \int_\Omega \rho_k \frac{|\bv-\bv_k|^2}{2h} \, dx
     + \int_\Omega 2 \eta(\p_k) |D\bv|^2 \, dx
     + \int_\Omega m(\p_k) |\nabla \mu|^2 \, dx \\
   & + \int_\Omega \widetilde{\Psi}_0'(A(\p)) \frac{A(\p)-A(\p_k)}{h} \, dx
     - \int_\Omega \widetilde{\kappa} \frac{A(\p)^2 - A(\p_k)^2}{2h} \, dx \\
   & + \frac{1}{h} \int_\Omega \nabla A(\p) \cdot \nabla (A(\p) - A(\p_k)) \, dx \\
   \geq & \int_\Omega \frac{\rho |\bv|^2 - \rho_k |\bv_k|^2}{2h} \, dx
     + \int_\Omega \rho_k \frac{|\bv-\bv_k|^2}{2h} \, dx
     + \int_\Omega 2 \eta(\p_k) |D\bv|^2 \, dx
     + \int_\Omega m(\p_k) |\nabla \mu|^2 \, dx \\ 
   & + \frac{1}{h} \int_\Omega \left( \widetilde{\Psi}_0(A(\p)) - \widetilde{\Psi}_0(A(\p_k)) \right) dx 
     - \int_\Omega \frac{\widetilde{\kappa}}{2} \, \frac{A(\p)^2 - A(\p_k)^2}{h} \, dx \\
   & + \frac{1}{h} \int_\Omega \frac{|\nabla A(\p) - \nabla A(\p_k)|^2}{2} \, dx
     + \frac{1}{h} \int_\Omega \left( \frac{|\nabla A(\p)|^2}{2} - \frac{|\nabla A(\p_k)|^2}{2} \right) dx \,,
\end{align*}
where we have used
\begin{align*}
 \widetilde{\Psi}_0'(A(\p)) \, (A(\p) - A(\p_k)) &\geq \widetilde{\Psi}_0(A(\p)) - \widetilde{\Psi}_0(A(\p_k)) \quad \mbox{and} \\
 \nabla A(\p) \cdot \nabla (A(\p) - A(\p_k)) &= \frac{|\nabla A(\p)|^2}{2} - \frac{|\nabla A(\p_k)|^2}{2} 
    + \frac{|\nabla A(\p) - \nabla A(\p_k)|^2}{2} \,.
\end{align*}
This leads to the claimed discrete energy estimate given by
\begin{align*} % \label{discreteenergyestimate} (number already defined above)
  E_{\mbox{\footnotesize tot}}&(\p,\bv) + \int_\Omega \rho_k \frac{|\bv - \bv_k|^2}{2} \, dx
    + \int_\Omega \frac{|\nabla A(\p) - \nabla A(\p_k)|^2}{2} \, dx \nonumber \\
    &+ h \int_\Omega 2 \eta(\p_k) |D\bv|^2 \, dx + h \int_\Omega m(\p_k) |\nabla \mu|^2 \, dx 
    \leq E_{\mbox{\footnotesize tot}}(\p_k,\bv_k) \,. 
 \end{align*}
In order to show existence of weak solutions we want to use the Leray-Schauder principle. 
With the abbreviation $H^2_n(\Omega) = \{u \in H^2(\Omega) \;|\; \partial_n u|_{\partial \Omega} = 0 \mbox{ on } \partial \Omega\}$ we define
operators $\La_{k},\F_k : X \to Y$, where
\begin{align*}
 X &= \left(H^1_0(\Omega)^d \cap L^2_\sigma(\Omega) \right) \times \D(\partial E) \times H^2_n(\Omega) \,, \\
 Y &= \left(H^1_0(\Omega)^d \cap L^2_\sigma(\Omega) \right)' \times L^2(\Omega) \times L^2(\Omega) \,.
\end{align*}
For $\bw = (\bv,\p,\mu) \in X$ we set
\begin{align*}
 \La_k (\bw) = \begin{pmatrix}
              L_k(\bv) \\
              -\mbox{div}(m(\p_k) \nabla \mu) + \int_\Omega \mu \, dx \\
              A(\p) + \partial \widetilde{E}(A(\p))
             \end{pmatrix} ,
\end{align*}
where 
\begin{align*}
 \left< L_k(\bv),\bpsi \right> &= \int_\Omega 2 \eta(\p_k) D\bv : D\bpsi \, dx \quad \mbox{for} \;
     \bpsi \in H^1_0(\Omega)^d \cap L^2_\sigma(\Omega) 
%\; \mbox{ and} \\
% \left<-\mbox{div}(m(\p_k) \nabla \mu) + \int_\Omega \mu , \eta \right> 
%    &= \int_\Omega m(\p_k) \nabla \mu \cdot \nabla \eta + \int_\Omega \mu \cdot \int_\Omega \eta \quad \mbox{for} \;
%     \eta \in H^1(\Omega) \,. %\\
% \left< A(\p) + \partial \widetilde{E}(A(\p)) , \eta \right> 
%    &= \int_\Omega A(\p) \, \eta 
%      + \int_\Omega \nabla A(\p) \cdot \nabla \eta
%      + \int_\Omega \widetilde{\Psi}_0'(A(\p)) \, \eta \quad \mbox{for} \; \eta \in H^1(\Omega) \,.
\end{align*}
and the second and third line are regarded pointwise.
Note that due to $\p \in \D(\partial E)$ it holds $A(\p) \in \D(\partial \widetilde{E})$ and therefore the last line 
in $\La_k(\bw)$ lies in $L^2(\Omega)$.
Furthermore for $\bw = (\bv,\p,\mu) \in X$ we define 
\begin{align*}
 \F_k(\bw) = \begin{pmatrix}
            - \frac{\rho \bv - \rho_k \bv_k}{h} - \mbox{div}(\rho_k \bv \otimes \bv) + \mu \nabla \p_k 
              - \left(\operatorname{div} \tbj - \frac{\rho - \rho_k}{h} - \bv \cdot \nabla \rho_k \right) \frac{\bv}{2}
              - \left( \tbj \cdot \nabla \right) \bv  \\
            -\frac{\p - \p_k}{h} - \bv \cdot \nabla \p_k + \int_\Omega \mu \, dx \rule{0cm}{0,5cm} \\
            A(\p) + \frac{\p - \p_k}{A(\p)-A(\p_k)} \,\mu + \tilde{\kappa} \frac{A(\p)+A(\p_k)}{2} \rule{0cm}{0,6cm}
           \end{pmatrix} .
\end{align*}
Then $\bw = (\bv,\p,\mu) \in X$ is a weak solution of the time discrete problem 
\eqref{timediscretizationline1}-\eqref{timediscretizationline3} if and only if
\begin{align*}
 \La_k (\bw) - \F_k(\bw) = 0 \,.
\end{align*}
Note that we used the equivalent version \eqref{equivtimediscretizationline1} instead
of \eqref{timediscretizationline1}, that is we substituted $\mbox{div}(\bv \otimes \tbj)$ accordingly. 

From standard theory of partial differential equations we get the invertibility of 
\[ L_k : H^1_0(\Omega)^d \cap L^2_\sigma(\Omega) \to \left(H^1_0(\Omega)^d \cap L^2_\sigma(\Omega)\right)' \,. \]
%and also the unique existence of a weak solution $\mu \in H^1(\Omega)$ of the second line, 
%which is afterwards shown to be in $H^2_n(\Omega)$ with the help of regularity theory.
%In fact this is done with the help of a bootstrapping argument, which we outline in what
%follows. 
Now we consider for given $f \in L^2(\Omega)$ the elliptic boundary value problem
\begin{align} \label{ell1}
 - \mbox{div}(m(\varphi_k) \nabla \mu) + \int_\Omega \mu \, dx = f \, \mbox{ in } \Omega \,, \quad
    \left. \partial_n \mu \right|_{\partial \Omega} = 0 \,.
\end{align}
With the Lemma of Lax-Milgram we get the unique existence of a weak solution $\mu \in H^1(\Omega)$ and 
we want to show with a bootstrapping argument that this solution fulfills $\mu \in H^2_n(\Omega)$.
At first, we observe that $\mu$ is also a weak solution to the following problem:
%This problem can also be written as 
\begin{align} \label{ell2}
 \Delta \mu = - \left( m(\varphi_k) \right)^{-1} \left( \nabla (m(\varphi_k)) \cdot \nabla \mu 
   + \int_\Omega \mu \, dx - f \right) 
   \, \mbox{ in } \Omega \,, \quad \left. \partial_n \mu \right|_{\partial \Omega} = 0 \,.
\end{align}
Since $\varphi_k \in H^2(\Omega)$ with $|\varphi_k| \leq 1$ almost everywhere in $\Omega$, we conclude
$\nabla \left( m(\varphi_k) \right) \in L^6(\Omega)$ and for the product of an $L^6$- with an 
$L^2$-function we get $\nabla \left( m(\varphi_k) \right) \cdot \nabla \mu \in L^{\frac{3}{2}}(\Omega)$. Using the bound
of $m(\varphi_k)$ from below by a positive constant, we get for some $\tilde{f} \in L^{\frac{3}{2}}(\Omega)$
that
\begin{align*}
 \Delta \mu = \tilde{f} \in L^{\frac{3}{2}}(\Omega) \,.
\end{align*}
Elliptic regularity theory now gives $\mu \in W^{2}_{\frac{3}{2}}(\Omega)$. In particular this means $\nabla \mu \in W^{1}_{\frac{3}{2}}(\Omega)$,
which embeds into $L^3(\Omega)$, so that we conclude for the product of an $L^6$- with an $L^3$-function
$\nabla \left( m(\varphi_k) \right) \cdot \nabla \mu \in L^2(\Omega)$. This shows 
\begin{align*}
 \Delta \mu = \tilde{f} \in L^{2}(\Omega) \,,
\end{align*}
and therefore finally $\mu \in H^2(\Omega)$ together with an estimate
\begin{align} \label{H2estimatemu}
 \|\mu\|_{H^2} \leq C \left( \|\mu\|_{H^1} + \|f\|_{L^2} \right) \,.
\end{align}
Since an $H^2$-solution of \eqref{ell2} also solves \eqref{ell1}, we obtain that the solution of 
\eqref{ell1} lies in $H^2_n(\Omega)$.

The operator $\partial \widetilde{E}$ is maximal monotone and hence
\begin{align*}
 I + \partial \widetilde{E} : \mathcal{D}(\partial \widetilde{E}) \rightarrow L^2(\Omega)
\end{align*}
is invertible. In addition we have continuity of the inverse considered as a mapping $L^2(\Omega) \to H^{2-s}(\Omega)$
for arbitrary $0 < s < \frac{1}{4}$. This is shown in the proof of Proposition 7.5.5 of Abels~\cite{Abe07} for a 
similar operator. For the convenience of the reader we give the details. Let $f_l \to f$ in $L^2(\Omega)$ for $l \to \infty$ with 
$f_l = u_l + \partial \widetilde{E}(u_l)$ and $f = u + \partial \widetilde{E}(u)$ be given. Then we have 
$u_l \to u$ in $H^1(\Omega)$ since
\begin{align*}
 \|u_l - u\|_{L^2}^2 + \|\nabla u_l - \nabla u\|_{L^2}^2 
   & \leq \|u_l - u\|_{L^2}^2 + \left( \partial \widetilde{E}(u_l) - \partial \widetilde{E}(u) , u_l - u \right)_{L^2} \\
   & \leq \|u_l + \partial \widetilde{E}(u_l) - (u + \partial \widetilde{E}(u)) \|_{L^2} \, \|u_l - u\|_{L^2} \\
   & \leq \frac{1}{2} \|f_l - f\|_{L^2}^2 + \frac{1}{2} \|u_l - u\|_{L^2}^2 \,.
\end{align*}
Moreover, because of the estimate \eqref{AWestimate}, $\|u_l\|_{H^2}$ is bounded and together with compactness we get 
the convergence $u_l \to u$ in $H^{2-s}(\Omega)$ for arbitrary $0<s<\frac{1}{4}$.

Since $A$ is invertible and we have the equivalence 
$A(\varphi) \in \mathcal{D}(\partial \widetilde{E})$ if and only if $\varphi
\in \mathcal{D}(\partial E)$, we also see that
\begin{align*}
 A(.) + \partial \widetilde{E}(A(.)) : \D(\partial E) \to L^2(\Omega)
\end{align*}
is invertible with continuous inverse considered from $L^2(\Omega) \to H^{2-s}(\Omega)$ as above. 
Altogether we obtain that $\La_k : X \to Y$ is invertible with inverse $\La_k^{-1} : Y \to X$.
Note that $X$ is not a Banach space since $\D(\partial E)$ includes inequality constraints. 

To get a continuous and even compact operator, we introduce for $0 < s < \frac{1}{4}$ the following Banach spaces
\begin{align*}
 \widetilde{X} &\ZuWeis \left( H_0^1(\Omega)^d \cap L_\sigma^2(\Omega) \right) 
                  \times H^{2 - s}(\Omega) \times H^2_n(\Omega) \,, \\
 \widetilde{Y} &\ZuWeis L^{\frac{3}{2}}(\Omega)^d \times W^{1}_{\frac{3}{2}}(\Omega) \times H^1(\Omega) \,.
\end{align*}
Then we obtain continuity of $\La_k^{-1} : Y \to \widetilde{X}$ from standard theory and with the above 
note concerning the continuity of the third line in $\La_k^{-1}$.
%To obtain continuity of $\La_k^{-1} : Y \to \widetilde{X}$, we refer to standard theory and to the work
%of Abels \cite{Abe07}, where it is shown in the proof of Proposition 7.5.5 for an analogue operator 
%that from $f_l \to f$ in $L^2(\Omega)$ the convergence $\p_l \to \p$ follows in $H^{2 - s}(\Omega)$.
%\textbf{(Rewrite the few lines nevertheless!)}

Finally due to $\widetilde{Y} \hookrightarrow \hookrightarrow Y$ the restriction 
$\La_k^{-1} : \widetilde{Y} \to \widetilde{X}$ is a compact operator. 

The next step is to observe that $\F_k : \widetilde{X} \to \widetilde{Y}$ is 
continuous and maps bounded sets into bounded sets. More precisely we obtain the following estimates (note that $\varphi_k$ and therefore $\rho_k$ 
lie in $H^2(\Omega)$):
\begin{align*}
 & \|\rho \bv\|_{L^{\frac{3}{2}}} \leq C \|\bv\|_{H^1} (\|\varphi\|_{L^2} + 1) \,, &
 & \|\mbox{div}(\rho_k \bv \otimes \bv) \|_{L^{\frac{3}{2}}}  \leq C_k \|\bv\|^2_{H^1} \,, \\
 & \|\mu \nabla \p_k \|_{L^{\frac{3}{2}}}  \leq C_k \|\mu\|_{L^2} \,, &
 & \|(\operatorname{div}\tbj) \bv\|_{L^{\frac{3}{2}}} \leq C_k \|\bv\|_{H^1} \|\mu\|_{H^2} \,, \\
 & \|(\tbj \cdot \nabla) \bv\|_{L^{\frac{3}{2}}} \leq C \|\bv\|_{H^1} \|\mu\|_{H^2} \,, &
 & \|\bv \cdot \nabla \p_k\|_{W^{1}_{\frac{3}{2}}} \leq C_k \|\bv\|_{H^1} \,, \\
 & \|A(\p)\|_{H^1} \leq C \|\p\|_{H^1} \,, &
 & \left\| \tfrac{\p - \p_k}{A(\p) - A(\p_k)} \mu \right\|_{H^1} \leq C_k \|\mu\|_{H^2} (\|\varphi\|_{H^1} + 1) \,. 
\end{align*}
%\textbf{(Some of the inequalities should be written down in detail, in part. $\mu \in H^2(\Omega) \hookrightarrow C^0(\overline{\Omega})$ 
%in the last inequality)} \\
In detail, these eight bounds are derived as follows:
\begin{enumerate}
 \item For $\rho \bv$, the estimate is clear due to $\bv \in H^1 \hookrightarrow L^6$. The ``+1'' appears,
       since $\rho$ depends affine linear on $\varphi$.
 \item For the estimate of $\operatorname{div}(\rho_k \bv \otimes \bv)$ in $L^{\frac{3}{2}}$, we have to estimate on the one hand 
       terms of the form $\rho_k \partial_l \bv_i \bv_j$
       in $L^{\frac{3}{2}}$, which are a product of functions in $L^\infty$, $L^2$ and $L^6$ and therefore are bounded in $L^{\frac{3}{2}}$.
       On the other hand we need terms of the form  $\partial_l \rho_k \bv_i \bv_j$ in $L^{\frac{3}{2}}$, which are a product of three functions
       in $L^6$ and therefore are bounded in particular also in $L^{\frac{3}{2}}$. 
 \item The estimate for $\mu \nabla \varphi_k$ in $L^{\frac{3}{2}}$ follows immediately from 
       $\mu \in L^6$ and $\partial_l \varphi_k \in L^2$.
 \item For the estimate of $(\operatorname{div}\tbj) \bv$ in $L^{\frac{3}{2}}$, where $\tbj = -\tfrac{\rho_2 - \rho_1}{2} m(\varphi_k) \nabla \mu$, 
       we have to estimate on the one hand terms of the form $m'(\varphi_k) \partial_i \varphi_k \partial_j \mu \bv_l$, which are a product
       of functions in $L^\infty$, $L^6$, $L^6$ and $L^6$ and therefore are bounded in particular in $L^{\frac{3}{2}}$. On the other hand
       we have to estimate terms of the form $m(\varphi_k) \partial_i \partial_j \mu \bv_l$, which are a product of functions
       in $L^\infty$, $L^2$ and $L^6$ and therefore are bounded in $L^{\frac{3}{2}}$.
 \item The estimate of $(\tbj \cdot \nabla) \bv$ in $L^{\frac{3}{2}}$ follows since the terms $m(\varphi_k) \partial_i \mu \partial_j \bv_l$ 
       are a product of functions in $L^\infty$, $L^6$ and $L^2$ and are therefore bounded in $L^{\frac{3}{2}}$.
 \item For the derivatives of $\bv \cdot \nabla \varphi_k$, we observe on the one hand that $\partial_i \bv_l$ is bounded in $L^2$ and 
       $\partial_j \varphi_k$ is bounded in $L^6$. On the other hand $\bv_i$ is bounded in $L^6$ and $\partial_i\partial_j \varphi_k$ is bounded 
       in $L^2$, so that in both cases the product is bounded in $L^{\frac{3}{2}}$.
 \item The estimate of $A(\varphi)$ follows from $|A(\varphi)| \leq K |\varphi|$ and $\partial_i (A(\varphi)) = A'(\varphi) \partial_i \varphi$ 
       almost everywhere, where $A'(\varphi) = \sqrt{a(\varphi)} \leq K$. 
 \item To estimate $\frac{\p - \p_k}{A(\p) - A(\p_k)} \mu$, we use first
       \begin{align*}
        \frac{\p - \p_k}{A(\p) - A(\p_k)} 
         = \left( \int_0^1 \underbrace{\sqrt{a(\tau \varphi + (1-\tau) \varphi_k)}}_{\geq c > 0} \, d\tau \right)^{-1} 
         \leq \frac{1}{c}
       \end{align*}
       and then from \eqref{helpestforquot} in the proof of Lemma \ref{derivativeforchempot} the result
       \begin{align*}
        \left| \nabla \left( \frac{A(\varphi) - A(\varphi_k)}{\varphi - \varphi_k} \right) \right| 
         \leq C(|\nabla \varphi| + 1) \,.
       \end{align*}
       In particular this gives a bound on 
       \begin{align*}
        \nabla \left( \frac{\p - \p_k}{A(\p) - A(\p_k)} \right) 
           = - \left(\frac{\p - \p_k}{A(\p) - A(\p_k)}\right)^2 
                 \nabla \left( \frac{A(\varphi) - A(\varphi_k)}{\varphi - \varphi_k} \right) .
       \end{align*}
       Together with $\mu \in H^2(\Omega) \hookrightarrow C^0(\overline{\Omega})$ we can bound 
       $\nabla \left( \frac{\p - \p_k}{A(\p) - A(\p_k)} \right) \mu$ in $L^2$ and of course also 
       $\frac{\p - \p_k}{A(\p) - A(\p_k)} \nabla \mu$ and $\frac{\p - \p_k}{A(\p) - A(\p_k)} \mu$ each in $L^2$.
\end{enumerate}
We want to use the Leray-Schauder principle on $\widetilde{Y}$ and rewrite therefore the 
identity $\La_k(\bw) - \F_k(\bw) = 0$ for a solution $\bw \in X$ of 
\eqref{timediscretizationline1}-\eqref{timediscretizationline3} into
\begin{align} \label{abstractproblem}
 \boldf - \F_k \circ \La_k^{-1}(\boldf) = 0 \quad \mbox{ for }\; \boldf = \La_k(\bw) \,.
\end{align}
We set $\K_k := \F_k \circ \La_k^{-1} : \widetilde{Y} \to \widetilde{Y}$ and note that $\K_k$ is 
a compact operator because $\La_k^{-1}$ is compact and $\F_k$ is continuous. 
Equation \eqref{abstractproblem} is then equivalent to finding a fixed point of $\K_k$, that is
\begin{align*}
 \boldf - \K_k(\boldf) = 0 \quad \Longleftrightarrow \quad \boldf = \K_k(\boldf) \,.
\end{align*}
To deduce the existence of such a fix point with the help of the Leray-Schauder 
principle (see for example Zeidler \cite{Zei92}), we have to show that
\begin{align} \label{requirementLeraySchauder}
 \exists R > 0 \, \mbox{ such that, if } \, \boldf \in \widetilde{Y} \, \mbox{ and } \, 0 \leq \lambda \leq 1 
  \, \mbox{ fulfills } \, \boldf = \lambda \K_k \boldf \,, \mbox{ then } \, \|\boldf \|_{\widetilde{Y}} \leq R \,.
\end{align}
So let $\boldf \in \widetilde{Y}$ and $0 \leq \lambda \leq 1$ with $\boldf = \lambda \K_k \boldf$. 
With $\bw = \La_k^{-1}(\boldf)$ we have
\begin{align*}
 \boldf = \lambda \K_k (\boldf) \quad \Longleftrightarrow \quad \La_k(\bw) - \lambda \F_k(\bw) = 0 \,,
\end{align*}
which is equivalent to the following weak formulation
\begin{align}
 &\int_\Omega 2 \eta(\p_k) D\bv : D\bpsi \, dx
  + \lambda \int_\Omega \frac{\rho \bv - \rho_k \bv_k}{h} \cdot \bpsi \, dx
  + \lambda \int_\Omega \mbox{div}(\rho_k \bv \otimes \bv) \cdot \bpsi \, dx \nonumber \\
  & \hspace*{10pt} + \lambda \int_\Omega \left( \mbox{div} \tbj - \frac{\rho - \rho_k}{h} - \bv \cdot \nabla \rho_k \right) \frac{\bv}{2} \cdot \bpsi \, dx
  + \lambda \int_\Omega \left( \tbj \cdot \nabla \right) \bv \cdot \bpsi \, dx
  = \lambda \int_\Omega \mu \nabla \p_k \cdot \bpsi \, dx \label{lambdaproblem1} \\
 & \mbox{for all } \; \bpsi \in H_0^1(\Omega)^d \cap L^2_\sigma(\Omega) \, \mbox{ and } \nonumber \\
 & \lambda \frac{\p - \p_k}{h} 
  + \lambda \bv \cdot \nabla \p_k  
  - \lambda \int_\Omega \mu \, dx
  = \mbox{div}(m(\p_k) \nabla \mu)  
  - \int_\Omega \mu \, dx \,, \label{lambdaproblem2} \\
 & A(\p) + \partial \widetilde{E}(A(\p)) 
  = \lambda A(\p) + \lambda \frac{\p - \p_k}{A(\p) - A(\p_k)} \mu 
                               + \lambda \widetilde{\kappa} \frac{A(\p) + A(\p_k)}{2} \,. \label{lambdaproblem3} 
\end{align}
First we derive an estimate for $\bw=(\bv,\p,\mu)$ in the norm of $\widetilde{X}$ and additionally an 
$L^2$-estimate of $\partial \widetilde{E}(A(\p))$, then we conclude the desired estimate for $\boldf$ 
due to the boundedness of $\F_k$. 

{\allowdisplaybreaks
Analogously as in the derivation of the energy estimate \eqref{discreteenergyestimate} we set 
$\bpsi = \bv$ in \eqref{lambdaproblem1}, test \eqref{lambdaproblem2} with $\mu$ and 
\eqref{lambdaproblem3} with $\frac{1}{h}(A(\p)-A(\p_k))$ to get with similar calculations
(we omit the integration element $dx$ for reasons of shortness):
\begin{align*}
 0 \; = \; & \lambda \frac{1}{h} \int_\Omega \left( \frac{\rho |\bv|^2}{2} - \frac{\rho_k |\bv_k|^2}{2} \right) 
      + \lambda \frac{1}{h} \int_\Omega \rho_k \frac{|\bv - \bv_k|^2}{2}
      + \int_\Omega 2 \eta(\p_k)  |D\bv|^2 
      + (1-\lambda) \left( \int_\Omega \mu \, dx \right)^2 \\
   &  + \int_\Omega m(\p_k) |\nabla \mu|^2 
      + (1-\lambda) \frac{1}{h} \int_\Omega A(\p) ( A(\p)-A(\p_k) )
      + \frac{1}{h} \int_\Omega \nabla A(\p) \cdot \left( \nabla A(\p) - \nabla A(\p_k) \right) \\
   &  + \frac{1}{h} \int_\Omega \widetilde{\Psi}_0'(A(\p)) (A(\p) - A(\p_k))
      - \lambda \frac{1}{h} \int_\Omega \widetilde{\kappa} \frac{A(\p)^2 - A(\p_k)^2}{2} \\
 \; \geq \; & \lambda \frac{1}{h} \int_\Omega \left( \frac{\rho |\bv|^2}{2} - \frac{\rho_k |\bv_k|^2}{2} \right) 
      + \lambda \frac{1}{h} \int_\Omega \rho_k \frac{|\bv - \bv_k|^2}{2}
      + \int_\Omega 2 \eta(\p_k)  |D\bv|^2 
      + (1-\lambda) \left( \int_\Omega \mu \, dx \right)^2 \\
   &  + \int_\Omega m(\p_k) |\nabla \mu|^2 
      + (1-\lambda) \frac{1}{h} \int_\Omega \left( \frac{A(\p)^2}{2} - \frac{A(\p_k)^2}{2} \right)
      + \frac{1}{h} \int_\Omega \left( \frac{|\nabla A(\p)|^2}{2} - \frac{|\nabla A(\p_k)|^2}{2} \right) \\
   &  + \frac{1}{h} \int_\Omega \left( \widetilde{\Psi}_0(A(\p)) - \widetilde{\Psi}_0(A(\p_k)) \right)
      - \lambda \frac{1}{h} \int_\Omega \widetilde{\kappa} \frac{A(\p)^2 - A(\p_k)^2}{2} \,.
\end{align*}}
This leads to the following estimate
\begin{align*}
 & h \int_\Omega 2 \eta(\p_k) |D\bv|^2 + h \int_\Omega m(\p_k) |\nabla \mu|^2 
 + \frac{1}{2} \int_\Omega |\nabla A(\p)|^2 + \int_\Omega \widetilde{\Psi}(A(\p)) + (1-\lambda)\left(\int_\Omega \mu \, dx \right)^2 \\
 & \leq \int_\Omega \frac{\rho_k |\bv_k|^2}{2} + \frac{1}{2} \int_\Omega A(\p_k)^2 
 + \frac{1}{2} \int_\Omega |\nabla A(\p_k)|^2 + \int_\Omega \widetilde{\Psi}_0(A(\p_k))
 + \int_\Omega |\widetilde{\kappa}| \frac{A(\p_k)^2}{2} \,,
\end{align*}
where we omitted the (nonnegative) terms $\lambda \int_\Omega \frac{\rho |\bv|^2}{2} \, dx$, 
$\lambda \int_\Omega \rho_k \frac{|\bv - \bv_k|^2}{2} \, dx$ and $(1-\lambda) \int_\Omega \frac{A(\p)^2}{2} \, dx$
on the left side, since due to the factor $\lambda$ resp. $(1-\lambda)$ they will not give 
a contribution to some estimate of $\|\bw\|_{\widetilde{X}}$ independent of $\lambda$. 
Note that due to $\bw=(\bv,\p,\mu) = \La_k^{-1}(\boldf) \in X$,
it holds that $\p \in \D(\partial E)$ and therefore $\p \in [-1,1]$ almost everywhere, which 
implies in particular $\rho \geq 0$. From this fact we also get boundedness of the term 
$\int_\Omega \widetilde{\Psi}(A(\p)) \, dx$, which can be estimated therefore on the right side. 
We also used the simple estimate $-\lambda \int_\Omega \widetilde{\kappa} \frac{A(\p_k)^2}{2} \, dx \leq
\lambda \int_\Omega |\widetilde{\kappa}| \frac{A(\p_k)^2}{2} \, dx$ and in addition we estimated every 
$\lambda$ resp. $(1-\lambda)$ on the right side against $1$.
%We remark that the term $-\lambda \int_\Omega \tilde{\kappa} \frac{A(\p_k)^2}{2}$ is nonnegative
%if we assume w.l.o.g. $\tilde{\kappa} \geq 0$, but this is not necessary. 

Using $|\nabla A(\p)|^2 = a(\p) |\nabla \p|^2$, we summarize the previous estimate to
\begin{align} \label{apriorilambda}
 (1-\lambda)\left(\int_\Omega \mu \, dx \right)^2 + h \int_\Omega 2 \eta(\p_k) |D\bv|^2 \, dx + h \int_\Omega m(\p_k) |\nabla \mu|^2 \, dx
 + \frac{1}{2} \int_\Omega a(\p) |\nabla \p|^2 \, dx \leq C_k \,. 
\end{align}
With the help of $\|\p\|_{L^\infty} \leq C$ due to $\p \in \D(\partial E)$, Korn's inequality for
$\bv \in H^1_0(\Omega)^d \cap L^2_\sigma(\Omega)$ and the fact that $\eta$, $m$ and $a$ are bounded
from below by a positive constant, we get the estimate
\begin{align} \label{apriorilambdashort}
 \sqrt{1-\lambda}\left|\int_\Omega \mu \, dx \right| + \|\bv\|_{H^1(\Omega)} + \|\nabla \mu\|_{L^2(\Omega)} + \|\p\|_{H^1(\Omega)} \leq C_k \,.
\end{align}
To get an estimate of the $L^2$-norm of the chemical potential $\mu$, we differ two cases. If 
$\lambda \in [\frac{1}{2},1]$, then we proceed with the simple estimate 
$\frac{1}{2} | \int_\Omega \mu \, dx| \leq \lambda |\int_\Omega \mu \, dx|$
and get as in the proof of Lemma \ref{derivativeforchempot} together with \eqref{apriorilambdashort} from 
equation \eqref{lambdaproblem3} the inequality
\begin{align*}
 \left| \int_\Omega \mu \, dx \right| \leq C_k \,.
\end{align*}
For $\lambda \in [0,\frac{1}{2})$ we use \eqref{apriorilambdashort} directly to get also here 
$\left| \int_\Omega \mu \, dx \right| \leq C_k$.
%For $\lambda \in [0,\frac{1}{2})$ we first integrate \eqref{lambdaproblem2} to get 
%$(\lambda - 1) |\Omega| \int_\Omega \mu \, dx = \lambda \int_\Omega \frac{\rho - \rho_k}{h} \, dx
%+ \lambda \int_\Omega \bv \cdot \nabla \varphi_k \, dx$. Then we use 
%$\frac{1}{2} |\int_\Omega \mu \, dx| \leq (1-\lambda) |\int_\Omega \mu \, dx|$ together with 
%\eqref{apriorilambdashort} to see also here $\left| \int_\Omega \mu \, dx \right| \leq C_k$. 
With this inequality for the mean value of the chemical 
potential $\mu$, we can improve estimate \eqref{apriorilambdashort} to
\begin{align} \label{apriorilambdashortfinal}
 \|\bv\|_{H^1(\Omega)} + \|\mu\|_{H^1(\Omega)} + \|\p\|_{H^1(\Omega)} \leq C_k \,.
\end{align}
Together with \eqref{H2estimatemu} we also get an estimate of the $H^2$-norm of the 
chemical potential $\mu$ given by
\begin{align} \label{apriorilambdashortfinally}
 \|\bv\|_{H^1(\Omega)} + \|\mu\|_{H^2(\Omega)} + \|\p\|_{H^1(\Omega)} \leq C_k \,.
\end{align}
Using identity \eqref{lambdaproblem3} pointwise we also get $\| \partial \widetilde{E}(A(\p)) \|_{L^2(\Omega)} \leq C_k$.
In summary this leads to
\begin{align*}
 \|\bw\|_{\widetilde{X}} + \|\partial \widetilde{E}(A(\p))\|_{L^2(\Omega)} 
 = \|(\bv,\p,\mu)\|_{\widetilde{X}} + \|\partial \widetilde{E}(A(\p))\|_{L^2(\Omega)} \leq C_k \,.
\end{align*}
To get finally an estimate of $\boldf = \La_k(\bw)$ in $\widetilde{Y}$ we use that
$\boldf - \lambda \F_k \La_k^{-1} (\boldf) = 0$ implies $\boldf = \lambda \F_k(\bw)$ and the fact that 
$\F_k : \widetilde{X} \to \widetilde{Y}$ maps bounded sets into bounded sets, which holds due to the 
above estimates for $\F_k$. This gives
\begin{align*}
 \|\boldf\|_{\widetilde{Y}} = \|\lambda \F_k(\bw) \|_{\widetilde{Y}} \leq C_k (\| \bw \|_{\widetilde{X}} + 1)
   \leq C_k \,,
\end{align*}
which was the remaining part to apply the Leray-Schauder principle as described above. Therefore 
we proved the existence of a weak solution of the time discrete problem 
\eqref{timediscretizationline1}-\eqref{timediscretizationline3}, which additionally fulfills the
discrete energy estimate \eqref{discreteenergyestimate}.
\end{proof}

\section{Proof of the Main Result}\label{sec:proof}
\subsection{Compactness in time}

To complete the proof of Theorem \ref{existenceweaksolution} we have to pass to the limit $h \to 0$ 
resp. $N \to \infty$ in our approximate solution. Therefore let $N \in \mathbb{N}$ be given and let 
$(\bv_{k+1},\p_{k+1},\mu_{k+1})$ be chosen successively as a solution of 
\eqref{timediscretizationline1}-\eqref{timediscretizationline3} with $h = \frac{1}{N}$ and $(\bv_0,\varphi_0^N)$ as initial value. 
Here we have to approximate the initial value $\varphi_0 \in H^1(\Omega)$ by $H^2$-functions $\varphi_0^N$
in order to apply Lemma \ref{timediscreteexistence} to the first step.
This is done with standard partial differential equation theory.
For example one could choose $u$ as the solution of the following heat equation
\begin{eqnarray*}
 \left\{ \begin{array}{rcll}
          \partial_t u - \Delta u &=& 0 & \mbox{in } \, \Omega \times (0,T) \,, \\
          u &=& \varphi_0 & \mbox{on } \, \Omega \times \{t=0\} \,, \\
          \left.\partial_\nu u\right|_{\partial \Omega} &=& 0 & \mbox{on } \, \partial \Omega \times (0,T) 
         \end{array}
 \right.
\end{eqnarray*}
and set $\varphi_0^N := \left.u\right|_{t=\frac{1}{N}}$ to get $\varphi_0^N \in H^2(\Omega)$, $| \varphi_0^N| \leq 1$ in $\Omega$
and $\varphi_0^N \to \varphi_0$ in $H^1(\Omega)$. 

Then we define $f^N(t)$ on $[-h,\infty)$ through $f^N(t) = f_k$ for $t \in [(k-1)h,kh)$, where
$k \in \mathbb{N}_0$ and $f \in \{\bv,\p,\mu\}$. 
Also set $\rho^N = \frac{1}{2}(\tilde{\rho}_1 + \tilde{\rho}_2) + \frac{1}{2}(\tilde{\rho}_2 - \tilde{\rho}_1) \varphi^N$.
In particular it holds $f^N((k-1)h) = f_k$, $f^N(kh) = f_{k+1}$ and $f^N(t) = f_{k+1}$ for $t \in [kh,(k+1)h)$. 
Additionally we define
\begin{align*}
 \left( \Delta^+_h f \right)(t) := f(t+h)-f(t) \,,\quad  & \left( \Delta^-_h f \right)(t) := f(t)-f(t-h) \,, \\
 \partial^+_{t,h} f(t) := \frac{1}{h} \left( \Delta_h^+ f \right)(t) \,, \quad &
 \partial^-_{t,h} f(t) := \frac{1}{h} \left( \Delta_h^- f \right)(t) \,, \\
 f_h := \left( \tau_h^\ast f\right)(t) = f(t-h) \,. \quad &
\end{align*}
Then for arbitrary $\bpsi \in \left( C^\infty_0(\Omega \times (0,\infty)) \right)^d$ with $\operatorname{div} \bpsi = 0$ we choose
$\widetilde{\bpsi} := \int_{kh}^{(k+1)h} \bpsi \, dt$ as test function in 
\eqref{timediscretizationline1} and sum over $k \in \mathbb{N}_0$ to get
\begin{align} 
  \int_0^\infty &\hspace*{-5pt} \int_\Omega \partial^-_{t,h}(\rho^N \bv^N) \cdot \bpsi \, dx \,dt
  + \int_0^\infty \hspace*{-5pt} \int_\Omega \mbox{div} \left(\rho^N_h \bv^N \otimes \bv^N \right) \cdot \bpsi \, dx \,dt
  + \int_0^\infty \hspace*{-5pt} \int_\Omega 2 \eta(\p^N_h) D\bv^N : D\bpsi \, dx \,dt \nonumber \\
 & - \int_0^\infty \hspace*{-5pt} \int_\Omega \left( \bv^N \otimes \tbj^N \right) : D\bpsi \, dx \,dt
   = \int_0^\infty \hspace*{-5pt} \int_\Omega \mu^N \nabla \p^N_h \cdot \bpsi \, dx \,dt \label{timeintegratedline1}
\end{align}
for all $\bpsi \in \left( C^\infty_0(\Omega \times (0,\infty)) \right)^d$ with $\operatorname{div} \bpsi = 0$.
The first term can be rewritten due to
\begin{align*}
 \int_0^\infty \hspace*{-5pt} \int_\Omega \partial^-_{t,h}(\rho^N \bv^N) \cdot \bpsi \, dx \,dt
   = - \int_0^\infty \hspace*{-5pt} \int_\Omega (\rho^N \bv^N) \cdot \partial^+_{t,h}\bpsi \, dx \,dt \,.
\end{align*}
Analogously we get 
\begin{align} \label{timeintegratedline2}
 \int_0^\infty \hspace*{-5pt} \int_\Omega \partial^-_{t,h} \p^N \, \zeta \, dx \,dt
  + \int_0^\infty \hspace*{-5pt} \int_{\Omega} \bv^N \p^N_h \cdot \nabla \zeta \, dx \,dt
  = \int_0^\infty \hspace*{-5pt} \int_\Omega m(\p^N_h) \nabla \mu^N \cdot \nabla \zeta \, dx \,dt
\end{align}
for all $\zeta \in C^\infty_0((0,\infty);C^1(\overline{\Omega}))$ and
\begin{align} 
 \frac{\Delta^-_h \p^N}{\Delta^-_h A(\p^N)} \, \mu^N 
  + \frac{\widetilde{\kappa}}{2} \left( A(\p^N) + A(\p^N_h) \right) 
 &= - \Delta A(\p^N) + \widetilde{\Psi}_0'(A(\p^N)) \label{timeintegratedline3}
\end{align}
holds pointwise in $\Omega \times (0,\infty)$ almost everywhere.

Let $E^N(t)$ be the piecewise linear interpolant of $E_{\mbox{\footnotesize tot}}(\p_k,\bv_k)$ at
$t_k = kh$ given by %for $t=kh$ it holds that $E^N(kh) = E_{\mbox{\footnotesize tot}}(\p_k,\bv_k)$ and
%for $t \in [kh,(k+1)h)$ we have
\begin{align*}
 E^N(t) = \frac{(k+1)h - t}{h} E_{\mbox{\footnotesize tot}}(\p_k,\bv_k) 
         + \frac{t - kh}{h} E_{\mbox{\footnotesize tot}}(\p_{k+1},\bv_{k+1}) \; \mbox{ for } \; t \in [kh,(k+1)h) \,.
\end{align*}
Also define for all $t \in (t_k,t_{k+1})$, $k \in \mathbb{N}_0$
\begin{align*}
 D^N(t) := \int_\Omega 2 \eta(\p_k) |D\bv_{k+1}|^2 \, dx + \int_\Omega m(\p_k) |\nabla \mu_{k+1}|^2 \, dx \,.
\end{align*}
Then the discrete energy estimate \eqref{discreteenergyestimate} implies
\begin{align} \label{inequalityforEandD}
 - \frac{d}{d t} E^N(t) = \frac{E_{\mbox{\footnotesize tot}}(\p_k,\bv_k) - E_{\mbox{\footnotesize tot}}(\p_{k+1},\bv_{k+1})}{h}
    \geq D^N(t) 
\end{align}
for all $t \in (t_k,t_{k+1})$, $k \in \mathbb{N}_0$. Multiplying this inequality by 
$\tau \in W^{1,1}(0,\infty)$ with $\tau \geq 0$, integrating and using integration by parts gives
\begin{align} \label{integratedinequEandD}
 E_{\mbox{\footnotesize tot}}(\p_0^N,\bv_0) \tau(0) + \int_0^\infty E^N(t) \, \tau'(t) \, dt \geq \int_0^\infty D^N(t) \, \tau(t) \, dt \,.
\end{align}
Integrating \eqref{inequalityforEandD} gives
\begin{align} %\label{inequforEN}
 E_{\mbox{\footnotesize tot}}(\p^N(t),\bv^N(t)) 
  &+ \int_s^t \int_\Omega \left( 2 \eta(\p^N_h) |D\bv^N|^2 + m(\p^N_h) |\nabla \mu^N|^2 \right) dx \,d\tau \nonumber \\
  &\leq E_{\mbox{\footnotesize tot}}(\p^N(s),\bv^N(s)) \label{inequforEN}
\end{align}
for all $0 \leq s \leq t < \infty$ with $s,t \in h\mathbb{N}_0$. 

Together with Lemma \ref{derivativeforchempot} and the fact that $E_{\mbox{\footnotesize tot}}(\varphi_0^N,\bv_0)$
is bounded this leads to the following bounds:
\begin{align} \label{timediscrbounds}
 \begin{array}{l}
  \bv^N \, \mbox{ is bounded in } \, L^2(0,\infty;H^1(\Omega)^d) \; \mbox{ and in } \; L^\infty(0,\infty; L^2(\Omega)^d) \,, \\
  \nabla \mu^N \, \mbox{ is bounded in } \, L^2(0,\infty;L^2(\Omega)^d) \,, \rule{0cm}{0,5cm}\\
  \varphi^N \, \mbox{ is bounded in } \, L^\infty(0,\infty;H^1(\Omega)) \, \mbox{ and } \rule{0cm}{0,5cm} \\
  \int_0^T \left| \int_\Omega \mu^N \, dx \right| dt \leq C(T) \, \mbox{ for all } \, 0<T<\infty \rule{0cm}{0,5cm}
 \end{array}
\end{align}
for a monotone function $C : \mathbb{R}^+ \to \mathbb{R}^+$. Using these bounds, we can pass to a subsequence to get
\begin{align*}
  \bv^N \rightharpoonup \bv \; &\mbox{ in } \, L^2(0,\infty;H^1(\Omega)^d) \,, \\
  \bv^N \rightharpoonup^\ast \bv \; &\mbox{ in } \, L^\infty(0,\infty;L^2(\Omega)^d) 
                               \cong \left( L^1(0,\infty;L^2(\Omega)^d) \right)' \,, \\
  \varphi^N \rightharpoonup^\ast \varphi \; &\mbox{ in } \, L^\infty(0,\infty;H^1(\Omega)) 
                               \cong \left( L^1(0,\infty;H^{1}(\Omega)) \right)' \,, \\
  \mu^N \rightharpoonup \mu \; &\mbox{ in } \, L^2(0,T;H^1(\Omega)) \, \mbox{ for all } \, 0<T<\infty \,, \\
  \nabla \mu^N \rightharpoonup \nabla \mu \; &\mbox{ in } \, L^2(0,\infty;L^2(\Omega)^d) \,.
\end{align*}
Here and in the following all limits are meant to be for suitable subsequences 
$N_k \to \infty$ (resp. $h_k \to 0$) for $k \to \infty$, unless otherwise stated. 

Now let $\widetilde{\varphi}^N$ be the piecewise linear interpolant of $\varphi^N(t_k)$, where 
$t_k = kh$, $h \in \mathbb{N}_0$, i.e. $\widetilde{\varphi}^N = \frac{1}{h} \chi_{[0,h]} \ast_t \varphi^N$,
where the convolution is only taken with respect to the time variable $t$. 
Then it holds that $\partial_t \widetilde{\varphi}^N = \partial_{t,h}^- \varphi^N$ and
\begin{align} \label{estofdiffint}
 \|\widetilde{\varphi}^N - \varphi^N \|_{H^{-1}(\Omega)} \leq h \|\partial_t \widetilde{\varphi}^N\|_{H^{-1}(\Omega)} \,.
\end{align}
From line \eqref{timeintegratedline2} we get that $\partial_t \widetilde{\varphi}^N \in L^2(0,\infty; H^{-1}(\Omega))$
is bounded, since $\bv^N \varphi^N$ and $\nabla \mu^N$ are both bounded in $L^2(0,\infty; L^2(\Omega)^d)$. 
%(Even $v^N \varphi^N \in L^2(0,\infty; L^2(\Omega))$ would be enough). 
Together with boundedness of $\widetilde{\varphi}^N$ in $L^\infty(0,\infty;H^1(\Omega))$, which 
follows from the boundedness of $\varphi^N$ in $L^\infty(0,\infty;H^1(\Omega))$, we get with the 
help of the lemma of Aubin-Lions \eqref{eq:AubinLions} the strong convergence 
\begin{align*}
 \widetilde{\varphi}^N \to \widetilde{\varphi} \; \mbox{ in } \; L^2(0,T;L^2(\Omega))
\end{align*}
for all $0<T<\infty$ for some $\widetilde{\varphi} \in L^\infty(0,\infty;L^2(\Omega))$. 
In particular it holds for a subsequence that $\widetilde{\varphi}^N \to \widetilde{\varphi}$
pointwise almost everywhere in $(0,\infty) \times \Omega$. 
Additionally by the above estimate \eqref{estofdiffint} we have 
\begin{align*}
 \widetilde{\varphi}^N - \varphi^N \to 0 \; \mbox{ in } \; L^2(0,\infty;H^{-1}(\Omega)) \,,
\end{align*}
which gives $\widetilde{\varphi} = \varphi$.
Furthermore, since $\widetilde{\varphi}^N \in H^1_{\mbox{\footnotesize uloc}}([0,\infty);H^{-1}(\Omega)) 
\cap L^2_{\mbox{\footnotesize uloc}}([0,\infty);H^1(\Omega))
\hookrightarrow BUC([0,\infty);L^2(\Omega))$ and $\widetilde{\varphi}^N \in L^\infty(0,\infty;H^1(\Omega))$
are bounded, it follows from Lemma \ref{lem:CwEmbedding} that $\varphi \in BC_w([0,\infty);H^1(\Omega))$.

%(\textbf{The spaces and the Lemma have to be stated precisely in the introduction.})

Finally it follows from the bound of $\widetilde{\varphi}^N$ in $H^1(0,T;H^{-1}(\Omega))$ and in $L^\infty(0,T;H^1(\Omega))$ 
for $0<T<\infty$ together with the fact that $\widetilde{\varphi}^N \to \varphi$ in $L^2(0,T;L^2(\Omega))$ that 
$\widetilde{\varphi}^N(0) \to \varphi(0)$ in $L^2(\Omega)$. But the left side equals $\varphi_0^N$, which converges
to $\varphi_0$ in $L^2(\Omega)$ so that we finally conclude $\varphi(0) = \varphi_0$.

Analogue observations can be done for $\rho^N$ since it depends affine linear on $\varphi^N$.

To show the convergence of \eqref{timeintegratedline3}, we observe that the right side is given 
by $\partial \widetilde{E}(A(\varphi^N))$ and this is bounded due to Lemma \ref{derivativeforchempot}
in $L^2(0,T;L^2(\Omega))$ for $0<T<\infty$. Moreover, the left side converges weakly in $L^2(0,T;L^2(\Omega))$ to
\begin{align*}
 f := a(\varphi)^{-\frac{1}{2}} \, \mu + \tilde{\kappa} A(\varphi) \,,
\end{align*}
which means that $\partial \widetilde{E}(A(\varphi^N)) \rightharpoonup f$ weakly in $L^2(0,T;L^2(\Omega))$.
If we now show that 
\begin{align} \label{limsupinequality}
 \limsup_{N \to \infty} \left< \partial \widetilde{E}(A(\varphi^N)) , A(\varphi^N) \right> \leq \left< f, A(\varphi) \right> \,, 
\end{align}
we can use the fact that $\partial \widetilde{E}$ is a maximal monotone operator and apply Prop. IV.1.6 in 
Showalter \cite{Sho97} to conclude that $\partial \widetilde{E}(A(\varphi)) = f$, which would finally lead to \eqref{weakline3}.

But the result \eqref{limsupinequality} does hold even with equality, since with the natural definition of $f^N$ 
as the left side of \eqref{timeintegratedline3}, we have
\begin{align*}
 \left< \partial \widetilde{E}(A(\varphi^N)) , A(\varphi^N) \right>
  = \left< f^N , A(\varphi^N) \right> \longrightarrow \left< f , A(\varphi) \right> 
\end{align*}
due to the strong convergence of $\varphi^N$ in $L^2(0,T;L^2(\Omega))$ and therefore also of $A(\varphi^N)$. 

Now we use the estimate \eqref{AWestimate} for $A(\varphi^N)$,
which together with Lemma \ref{derivativeforchempot} gives boundedness of $A(\varphi^N)$ in 
$L^2(0,T;H^2(\Omega))$. Then the same holds true for $\varphi^N$ and due to the strong convergence of 
$\varphi^N$ in $L^2(0,T;H^{-1}(\Omega))$ we conclude with an interpolation argument even the strong 
convergence 
\begin{align} \label{strongphiinL2H1}
 \varphi^N \to \varphi \; \mbox{ in } \;  L^2(0,T;H^1(\Omega)) \,.
\end{align}
The next step is to show strong convergence $\bv^N \to \bv$ in $L^2(0,T;L^2(\Omega)^d)$ for all
$0<T<\infty$ to conclude a convergence pointwise almost everywhere.
As above let $\widetilde{\rho \bv}^N$ be the piecewise linear interpolant of $\left(\rho^N \bv^N\right)(t_k)$, where 
$t_k = kh$, $h \in \mathbb{N}_0$. Then it holds that 
$\partial_t \left(\widetilde{\rho \bv}^N\right) = \partial_{t,h}^- \left(\rho^N \bv^N\right)$.% and
%\begin{align} \label{estofdiffint2}
% \|\widetilde{\rho \bv}^N - \rho^N \bv^N \|_{H^{-1}(\Omega)} 
%    \leq h \| \partial_t \left(\widetilde{\rho \bv}^N\right) \|_{H^{-1}(\Omega)} \,.
%\end{align}

With the help of the projection $\mathbb{P}_\sigma$ onto $L^2_\sigma(\Omega)$ we get from line \eqref{timeintegratedline1} that 
$\partial_t \left( \mathbb{P}_\sigma(\widetilde{\rho \bv}^N) \right)$ is bounded in 
$L^{\frac{8}{7}}(0,T;W^{-1}_4(\Omega))$ for $0<T<\infty$, since due to some known interpolation inequalities 
we have the following bounds 
\begin{align*}
 \rho^N_h \bv^N \otimes \bv^N & \; \mbox{ is bounded in } \; L^2(0,T;L^\frac{3}{2}(\Omega)) \,, \\
 D\bv^N                       & \; \mbox{ is bounded in } \; L^{2}(0,T;L^2(\Omega)) \,, \\
 \bv^N \otimes \nabla \mu^N          & \; \mbox{ is bounded in } \; L^{\frac{8}{7}}(0,T;L^{\frac{4}{3}}(\Omega)) \,, \\
 \mu^N \nabla \varphi^N_h     & \; \mbox{ is bounded in } \; L^2(0,T;L^\frac{3}{2}(\Omega)) \,.
\end{align*}
%\textbf{More details, which will be probably deleted in the preprint:} 
%%%%%%%%%%%%%%%%%%%%%%%%%%%%%%%%%%%%%%%%%%%%%%%%%%%%%%%%%%%%%%%%%%%%%%%%%%%%
%%%%%%%%%%%%%%%%%%%%%%%%%%%%%%%%%%%%%%%%%%%%%%%%%%%%%%%%%%%%%%%%%%%%%%%%%%%%
%%%%%%%%%%%%%%%%%%%%%%%%%%%%%%%%%%%%%%%%%%%%%%%%%%%%%%%%%%%%%%%%%%%%%%%%%%%%
%%%%%%%%%%%%%%%%%%%%%%%%%%%%%%%%%%%%%%%%%%%%%%%%%%%%%%%%%%%%%%%%%%%%%%%%%%%%
In detail, these bounds are derived as follows, where for abbreviation ``$\in$'' always means bounded independent of $N$ in the 
corresponding space and we omit the time and space variables. 
\begin{enumerate}
 \item Due to $\bv^N \in L^\infty(L^2) \cap L^2(L^6)$, we get for products
       $\bv^N_i \bv^N_j \in L^2(L^{\frac{3}{2}})$ and from $\rho^N_h \in L^\infty(\Omega_T)$ this also holds 
       for $\rho^N_h \bv^N \otimes \bv^N$.
 \item From $\bv \in L^2(H^1)$ we immediately get $D\bv^N \in L^2(L^2)$. 
 \item Due to $\bv^N \in L^\infty(L^2) \cap L^2(L^6)$ and $\nabla \mu^N \in L^2(L^2)$ we get $\bv^N_i \partial_j \mu^N \in L^2(L^1)
       \cap L^1(L^{\frac{3}{2}})$. Now, if $I \subset \mathbb{R}$ is an interval,
       $(X_0,X_1)$ is an interpolation couple of Banach spaces, 
       and $X = (X_0,X_1)_{\theta, q}$ is the real interpolation space of type $\theta \in (0,1)$ with 
       $1 \leq q \leq \infty$, we have the embedding
       \begin{align} \label{interpolembedding}
        L^{p_0}(I;X_0) \cap L^{p_1}(I;X_1) \hookrightarrow L^p(I;X) \,, 
         \quad  \mbox{where } \, \frac{1}{p} = \frac{1-\theta}{p_0} + \frac{\theta}{p_1} \,,
       \end{align}
       where $1 \leq p_0,p_1 \leq \infty$.
       In addition, we use the description of real interpolation spaces with $0<\theta < 1$ 
       for Lebesgue-spaces with $1 \leq q_0,q_1 \leq \infty$ given trough
       \begin{align} \label{realinterpollebesgue}
         (L^{q_0}(\Omega),L^{q_1}(\Omega))_{\theta,q} = L^q(\Omega) \,, 
          \quad  \mbox{where } \, \frac{1}{q} = \frac{1-\theta}{q_0} + \frac{\theta}{q_1} \,.
       \end{align}
       For $X_0 = L^1(\Omega)$ and $X_1 = L^{\frac{3}{2}}(\Omega)$ we get with $\theta = \frac{3}{4}$ that
       $q = \frac{4}{3}$ and therefore $X = (X_0,X_1)_{\theta,q} = L^\frac{4}{3}(\Omega)$. 
       Furthermore with $p_0 = 2$ and $p_1 = 1$ we get $p = \frac{8}{7}$, which leads to the bound
       $\bv^N_i \partial_j \mu^N \in L^{\frac{8}{7}}(L^{\frac{4}{3}})$.
 \item From $\mu^N \in L^2(L^6)$ and $\nabla \varphi^N_h \in L^\infty(L^2)$ we get finally 
       $\mu^N \nabla \varphi^N_h \in L^2(L^{\frac{3}{2}})$. 
\end{enumerate}
This means in particular that all four terms are bounded in $L^{\frac{8}{7}}(0,T;L^{\frac{4}{3}}(\Omega))$.
Therefore we can allow in \eqref{timeintegratedline1} for test functions with
$\bpsi$, $\nabla \bpsi \in \left(L^\frac{8}{7}(L^{\frac{4}{3}})\right)' = L^8(L^4)$, i.e. $\bpsi \in 
L^8(W^{1}_{4})$. This implies the bound for $\partial_t \left( \mathbb{P}_\sigma(\widetilde{\rho \bv}^N) \right)$
in $\left(L^8(W^{1}_{4})\right)' = L^\frac{8}{7}(W^{-1}_{4})$.

Additionally $\mathbb{P}_\sigma (\widetilde{\rho \bv}^N) \in L^2(0,T;H^1(\Omega)^d)$ is bounded 
and we can therefore conclude with the help of the Lemma of Aubin-Lions \eqref{eq:AubinLions} the strong convergence 
\begin{align*}
 \mathbb{P}_\sigma (\widetilde{\rho \bv}^N) \to \bw
  \; \mbox{ in } \; L^2(0,T;L^2(\Omega)^d)
\end{align*}
for all $0<T<\infty$ for some $\bw \in L^\infty(0,\infty;L^2(\Omega)^d)$. 

Since the projection $\mathbb{P}_\sigma : L^2(0,T;L^2(\Omega)^d) \to L^2(0,T;L^2_\sigma(\Omega))$ is weakly continuous, we 
conclude from the weak convergence $\widetilde{\rho \bv}^N \rightharpoonup \rho \bv$ 
in $L^2(0,T;L^2(\Omega))$ that $\bw = \mathbb{P}_\sigma (\rho \bv)$.
%Additionally by the above estimate \eqref{estofdiffint2} we have 
%\begin{align*}
% \mathbb{P}_\sigma \left( \widetilde{\rho \bv}^N - \rho^N \bv^N \right) \to 0 \; \mbox{ in } \; L^2(0,T;H^{-1}(\Omega)) \,,
%\end{align*}
%which gives $\widetilde{\rho \bv} = \rho \bv$. 

Now we derive the strong convergence 
$\bv^N \to \bv$ in $L^2(0,T;L^2(\Omega)^d)$ with the help of the following observations:
\begin{align*}
 \int_0^T \int_\Omega \rho^N |\bv^N|^2 = \int_0^T \int_\Omega \mathbb{P}_\sigma (\rho^N \bv^N ) \cdot \bv^N
 \longrightarrow \int_0^T \int_\Omega \mathbb{P}_\sigma(\rho \bv) \cdot \bv 
  = \int_0^T \int_\Omega \rho \, |\bv|^2 \,,
\end{align*}
where we used the strong convergence of $\mathbb{P}_\sigma(\rho^N \bv^N)$ in $L^2(0,T;L^2(\Omega)^d)$ and
the weak convergence of $\bv^N$ in $L^2(0,T;L^2(\Omega)^d)$.
This gives $(\rho^N)^{\frac{1}{2}} \bv^N \to (\rho)^{\frac{1}{2}} \bv$ in $L^2(0,T;L^2(\Omega)^d)$
and from above we know that 
\begin{align*}
 \rho^N \to \rho \; \mbox{ almost everywhere in } \, (0,\infty) \times \Omega \; \mbox{ and } \;
  |\rho^N| \geq c > 0 \,,
\end{align*}
so that we can conclude
\begin{align*}
 \bv^N = (\rho^N)^{-\frac{1}{2}} \left( (\rho^N)^{\frac{1}{2}} \bv^N \right)
  \to \bv \; \mbox{ in } \; L^2(0,T;L^2(\Omega)^d) \,.
\end{align*}
This means in particular that $\bv^N \to \bv$ pointwise almost everywhere in $(0,\infty) \times \Omega$
(for a subsequence). We note that similar ideas to prove strong convergence of $\bv^N$ have been used 
by P.-L.~Lions~\cite[Section 2.1]{Lio96}.

Using these convergence results together with the fact that for all divergence free $\bpsi$ the 
following convergence holds
\begin{align*}
 \int_0^T \int_\Omega \mu^N \nabla \varphi^N_h \cdot \bpsi = -\int_0^T \int_\Omega \nabla \mu^N \varphi^N_h \cdot \bpsi \longrightarrow 
  -\int_0^T \int_\Omega \nabla \mu \varphi \cdot \bpsi = \int_0^T \int_\Omega \mu \nabla \varphi \cdot \bpsi \,,
\end{align*}
we can pass to the limit in the equations \eqref{timeintegratedline1},
\eqref{timeintegratedline2} to get \eqref{weakline1}, \eqref{weakline2}.

\subsection{Initial data for the velocity $\bv$}

\begin{lemma} \label{freedivlemma}
 Let $\bv$, $\tilde{\bv} \in L^2_\sigma(\Omega) \cap H^1_0(\Omega)^d$ 
 and $\rho \in L^\infty(\Omega)$ with $\rho \geq c > 0$ such that
 \begin{eqnarray*}
  \int_{\Omega} \rho \bv \cdot  \bpsi \, dx = \int_{\Omega} \rho \tilde{\bv} \cdot \bpsi \, dx \quad \mbox{ for all } \; 
   \bpsi \in C_{0,\sigma}^\infty(\Omega) \,.
 \end{eqnarray*}
 Then it holds that $\bv = \tilde{\bv}$ almost everywhere in $\Omega$. 
\end{lemma}
\begin{proof}
 By approximation the identity also holds for $\bpsi := \bv - \tilde{\bv} \in L^2_\sigma(\Omega) \cap H^1_0(\Omega)$
and we get $\int_\Omega \rho_0 |\bv - \tilde{\bv}|^2 \, dx = 0$, which due to the boundedness of $\rho_0$ from below
by a positive constant gives finally $\bv = \tilde{\bv}$ in $L^2(\Omega)^d$.
\end{proof}

\iffalse \begin{proof} By assumption we can write $\rho (\bv - \tilde{\bv}) = \nabla p$ for some $p \in H^1(\Omega)$. Since
$\rho$ is bounded from below by a positive constant, this is equivalent to $(\bv - \tilde{\bv}) = \frac{1}{\rho} \nabla p$.
Then $p$ is a weak solution of the problem
\begin{eqnarray*}
 \left\{ \begin{array}{rcll}
          \operatorname{div}(\frac{1}{\rho} \nabla p) &=& 0 & \mbox{in } \, \Omega \,, \\
          \nabla p \cdot n &=& 0 & \mbox{on } \, \partial \Omega \,,
         \end{array}
 \right.
\end{eqnarray*}
which means 
\begin{align*}
 \int_\Omega \frac{1}{\rho} \nabla p \cdot \nabla \varphi \, dx = 0 \; \mbox{ for all } \; \varphi \in H^1(\Omega) \,. 
\end{align*}
Since its solution is just $p = const$, we get $\rho (\bv - \tilde{\bv}) = 0$ and finally $\bv = \tilde{\bv}$
and thus we proved the Lemma. 
\end{proof} \fi

%\vskip 5pt 
With the following arguments we show that the velocity $\bv$ fulfills the initial condition
\eqref{weakline4}. At first we derive weak continuity in time for the projection 
of $\rho \bv$ onto the divergence-free vector fields. Therefore let $\bw^N := \mathbb{P}_\sigma(\widetilde{\rho \bv}^N)$
and from the arguments above we deduce the boundedness of
\begin{align*}
 & \bw^N \; \mbox{ in } \; W^{1}_{\frac{8}{7},\mbox{\scriptsize uloc}}([0,\infty);W^{-1}_{4}(\Omega))
   \hookrightarrow BUC([0,\infty);W^{-1}_{4}(\Omega)) \; \mbox{ and } \\
 & \bw^N \; \mbox{ in } \; L^\infty(0,\infty;L^2(\Omega)^d)  \,.
\end{align*}
Due to these bounds, we concluded above that $\bw^N \to \bw$ in $L^2(0,T;L^2(\Omega)^d)$ and with the help of
Lemma \ref{lem:CwEmbedding} we deduce for the limit that $\bw \in BC_w([0,\infty);L^2(\Omega)^d)$. 
Now let $0<T<\infty$ and consider the auxiliary problem for a function $p \in L^2(0,T;H^1_{(0)}(\Omega))$
given by
\begin{equation*} 
 \left\{ \begin{array}{rcll}
          -\operatorname{div}( \frac{1}{\rho(t)} \nabla p(t) ) &=& \operatorname{div} (\frac{1}{\rho(t)} \bw(t)) & \mbox{in } \; \Omega \,, \\
          \nabla p(t) \cdot n &=& 0 & \mbox{on } \; \partial \Omega \,.
         \end{array}
 \right. \tag{$Aux$}
\end{equation*}
The corresponding weak version reads as
\begin{eqnarray}
 -\int_\Omega \frac{1}{\rho(t)} \nabla p(t) \cdot \nabla \tau \, dx
  = \int_\Omega \frac{1}{\rho(t)} \bw(t) \cdot \nabla \tau \, dx
   \quad \mbox{ for all } \; \tau \in H^1(\Omega)  \label{Auxweak}
\end{eqnarray}
for almost all $t$ in $(0,T)$. By the Lemma of Lax-Milgram there exists a unique solution and it fulfills the estimate
\begin{eqnarray}
 \| \nabla p(t) \|_{L^2(\Omega)} \leq C \| \bw(t) \|_{L^2(\Omega)} \,, \label{Auxestimate}
\end{eqnarray}
where $C$ is independent of $t$.
Here we used the fact that $\rho \in L^\infty(\Omega_T)$ is bounded from below by a positive constant.

\noindent
Claim: $\nabla p \in BC_w([0,T];L^2(\Omega)^d)$. 

\noindent
To see this claim, let $t_n$, $t \in [0,T]$ for $n \in \mathbb{N}$ with $t_n \to t$ be given. From 
$\bw \in BC_w([0,\infty);L^2(\Omega)^d)$ we conclude that in particular $\bw(t_n)$ is bounded in $L^2$ and
from the last estimate \eqref{Auxestimate} also $\nabla p(t_n)$ is bounded in $L^2$. 
Therefore we get the weak convergence $\nabla p(t_{n_k}) \rightharpoonup \nabla q \in L^2(\Omega)^d$ for $k \to \infty$
at least for a subsequence, since $\{\nabla q \;|\; q \in H^1_{(0)}(\Omega)\}$ is a closed subspace in $L^2(\Omega)^d$. 
But due to the unique solvability of \eqref{Auxestimate} we arrive at 
$ \nabla q = \nabla p(t)$. Since this argumentation holds for any weakly convergent subsequence, we even have 
$\nabla p(t_n) \rightharpoonup \nabla p(t)$. 

In addition we already know due to the weak 
continuity of $\mathbb{P}_\sigma$ and due to the weak convergence $\widetilde{\rho \bv}^N \to \rho \bv$ in 
$L^2(0,T;L^2(\Omega)^d)$ for every $0<T<\infty$, that $\bw = \mathbb{P}_\sigma (\rho \bv)$.
In particular this gives
\begin{eqnarray*}
 \bw(t) = \mathbb{P}_\sigma\left( \rho(t) \bv(t) \right) \;\; \mbox{ almost everywhere in } (0,\infty) \,.
\end{eqnarray*}
By characterization \eqref{eq:WeakHelmholtz} of the projection $\mathbb{P}_\sigma$ onto $L^2_\sigma(\Omega)$, the last 
identity means 
\begin{eqnarray} \label{lineforbw}
 \bw(t) = \rho(t) \bv(t) - \nabla \tilde{p}(t) \;\; \mbox{ almost everywhere in } (0,\infty) \,,
\end{eqnarray}
where $\tilde{p}(t)$ is a weak solution of 
\begin{eqnarray*}
 \left\{ \begin{array}{rcll}
          \operatorname{div}(\nabla \tilde{p}(t)) &=& \operatorname{div}(\rho(t) w(t)) & \mbox{in } \; \Omega \,, \\
          \nabla \tilde{p}(t) \cdot n &=& 0 & \mbox{on } \; \partial \Omega \,.
         \end{array}
 \right.
\end{eqnarray*}
Dividing \eqref{lineforbw} by $\rho(t)$, we get 
\begin{eqnarray} \label{lineforbv}
 \bv(t) = \frac{1}{\rho(t)} \bw(t) + \frac{1}{\rho(t)} \nabla \tilde{p}(t) \;\; \mbox{ almost everywhere in } (0,\infty) \,,
\end{eqnarray}
which after multiplication with $\nabla \tau$ for $\tau \in H^1(\Omega)$ and integration gives
\begin{eqnarray*}
 \int_\Omega \frac{1}{\rho(t)} \bw(t) \cdot \nabla \tau \, dx 
   = - \int_\Omega \frac{1}{\rho(t)} \nabla \tilde{p}(t) \cdot \nabla \tau \, dx 
     \quad \mbox{ for all } \; \tau \in H^1(\Omega) 
\end{eqnarray*}
for almost all $t \in (0,\infty)$. But this is exactly the weak version \eqref{Auxweak} of the above auxiliary problem $(Aux)$ and 
we get $\nabla \tilde{p}(t) = \nabla p(t)$ for almost every $t$ due to the unique solvability.
By a redefinition of $\nabla \tilde{p}$ on a set of measure zero we get the continuity
$\nabla \tilde{p} \in BC_w([0,T];L^2(\Omega)^d)$, since $\nabla p$ has this property. 
Again with a redefinition of $\bv$ through equation \eqref{lineforbv} on a set of measure zero, 
we get finally $\bv \in BC_w([0,T];L^2(\Omega)^d)$, where we used the fact $\rho \in BUC([0,T];L^2(\Omega))$
with $\rho \geq c > 0$ is already known from above. 

It remains to show that $\bv$ attains the initial value $\bv_0$.  Due to the above bound of 
$\mathbb{P}_\sigma(\widetilde{\rho \bv}^N)$ in $W^{1}_{\frac{8}{7}}(0,T;W^{-1}_{4}(\Omega))
\hookrightarrow BUC([0,T];W^{-1}_{4}(\Omega))$ we get for arbitrary $\bpsi \in C^\infty_{0,\sigma}(\Omega)$ 
the convergence 
\begin{eqnarray*}
 \int_\Omega \mathbb{P}_\sigma(\widetilde{\rho \bv}^N)(0) \cdot \bpsi \, dx 
 \longrightarrow \int_\Omega \mathbb{P}_\sigma(\rho_0 \bv(0)) \cdot \bpsi
  = \int_\Omega \rho_0 \bv(0) \cdot \bpsi \, dx \,.
\end{eqnarray*}
By definition of the time-discrete functions we also have $\mathbb{P}_\sigma(\widetilde{\rho \bv}^N)|_{t=0}
= \mathbb{P}_\sigma (\rho_0 \bv _0)$, which together with the last convergence yields
\begin{eqnarray*}
 \int_\Omega \rho_0 \bv_0 \cdot \bpsi \, dx = \int_\Omega \rho_0 \bv(0) \cdot \bpsi \, dx 
  \quad \mbox{ for all } \; \bpsi \in C^\infty_{0,\sigma}(\Omega) \,.
\end{eqnarray*}
From Lemma \ref{freedivlemma} we get finally $\bv(0) = \bv_0$ in $L^2(\Omega)^d$.

\subsection{Energy inequality}

Finally we can finish the proof by showing the energy inequality \eqref{weakline5}.
Since $\bv^N(t) \to \bv(t)$ in $L^2(\Omega)^d$ and $\varphi^N(t) \to \varphi(t)$ in $H^1(\Omega)$ for
almost every $t \in (0,\infty)$ (for a subsequence), which follows from the strong convergences of 
$\bv^N$ and $\varphi^N$, it holds that
\begin{align*}
 E^N(t) \to E_{\mbox{\footnotesize tot}}(\varphi(t),\bv(t)) \; \mbox{ for almost all } \, t \in (0,\infty) \,.
\end{align*}
Moreover, by lower semicontinuity of norms and almost everywhere convergence of $\varphi^N$ to $\varphi$,
the inequality 
\begin{align*}
 \liminf_{N \to \infty} \int_0^\infty D^N(t) \tau(t) \, dt \geq \int_0^\infty D(t) \tau(t) \, dt 
\end{align*}
for all $\tau \in W^{1,1}(0,\infty)$ with $\tau \geq 0$ holds, where
\begin{align*}
 D(t) := \int_\Omega 2 \eta(\varphi) |D\bv|^2 \, dx + \int_\Omega m(\varphi) |\nabla \mu|^2 \, dx \,.
\end{align*}
Hence, passing to the limit in \eqref{integratedinequEandD}, we obtain
\begin{align} \label{lastinequalityfinally}
 E_{\mbox{\footnotesize tot}}(\varphi_0,\bv_0) \tau(0) 
   + \int_0^\infty E_{\mbox{\footnotesize tot}}(\varphi(t),\bv(t)) \, \tau(t) \, dt 
  \geq \int_0^\infty D(t) \tau(t) \, dt 
\end{align}
for all $\tau \in W^{1,1}(0,\infty)$ with $\tau \geq 0$. With the help of Lemma \ref{lem:EnergyEstim} 
we obtain the energy estimate~\eqref{weakline5}. 
\hfill $\Box$

%%%%%%%%%%%%%%%%%%%%%%%%%%%%%%%%%%%%%%%%%%%%%%%%%%%%%%%%%%%%%%%%%%%%%%%%%%%%%%%%%%%%%%%%%%%%%%%%%%%%%%%%%%%%%%%%%%%%%%%%%
%%%%%%%%%%%%%%%%%%%%%%%%%%%%%%%%%%%%%%%%%%%%%%%%%%%%%%%%%%%%%%%%%%%%%%%%%%%%%%%%%%%%%%%%%%%%%%%%%%%%%%%%%%%%%%%%%%%%%%%%%
%%%%%%%%%%%%%%%%%%%%%%%%%%%%%%%%%%%%%%%%%%%%%%%%%%%%%%%%%%%%%%%%%%%%%%%%%%%%%%%%%%%%%%%%%%%%%%%%%%%%%%%%%%%%%%%%%%%%%%%%%
%%%%%%%%%%%%%%%%%%%%%%%%%%%%%%%%%%%%%%%%%%%%%%%%%%%%%%%%%%%%%%%%%%%%%%%%%%%%%%%%%%%%%%%%%%%%%%%%%%%%%%%%%%%%%%%%%%%%%%%%%
%%%%%%%%%%%%%%%%%%%%%%%%%%%%%%%%%%%%%%%%%%%%%%%%%%%%%%%%%%%%%%%%%%%%%%%%%%%%%%%%%%%%%%%%%%%%%%%%%%%%%%%%%%%%%%%%%%%%%%%%%
%%%%%%%%%%%%%%%%%%%%%%%%%%%%%%%%%%%%%%%%%%%%%%%%%%%%%%%%%%%%%%%%%%%%%%%%%%%%%%%%%%%%%%%%%%%%%%%%%%%%%%%%%%%%%%%%%%%%%%%%%

\renewcommand{\theequation}{A.\arabic{equation}}

\appendix
\section{Appendix}
In this last section, we want to discuss existence of weak solutions for a related model discussed in Remark 2.2 of
Abels, Garcke, Gr\"un \cite{AGG11}, which is derived in detail in \cite{AGG10}. 
Although this model is not frame indifferent in the standard manner, it can be used to approximate the sharp interface
model, see \cite{AGG10}. In addition it has the mathematical 
interesting term $|\bv|^2$ in the equation for the chemical potential and we want to give a short description how the methods from this work 
can be adapted to derive an existence result also in that case. After reformulation of the pressure term with the 
previous notation, this problem is given as
\begin{align}  
  \partial_t(\rho \bv) + \operatorname{div} \left( \rho \bv \otimes \bv \right) - \operatorname{div} \left( 2\eta(\p) D\bv \right) + \nabla g 
     &= \mu \nabla \p + \frac{|\bv|^2}{2} \nabla \rho  & \mbox{in } \; Q , \label{old1} \\
  \operatorname{div} \bv &= 0 & \mbox{in } \; Q , \label{old2} \\
  \partial_t \varphi + \bv \cdot \nabla \varphi &= \operatorname{div} \left(m(\varphi) \nabla \mu\right) & \mbox{in } \; Q , \label{old3} \\
  a^{-\frac{1}{2}}(\p) \left( \mu + \frac{\partial \rho}{\partial \varphi} \frac{|\bv|^2}{2} \right) + \widetilde{\kappa} A(\p)  
   &= \widetilde{\Psi}_0'(A(\p)) - \Delta A(\p)  & \mbox{in } \; Q , \label{old4} \\
  \bv |_{\partial \Omega} & = 0 & \mbox{on } \, S , \label{old5} \\ 
 \partial_n \varphi |_{\partial \Omega} = \partial_n \mu |_{\partial \Omega} & = 0 & \mbox{on } \, S , \label{old6} \\ 
 (\bv,\varphi) |_{t=0} & = (\bv_0,\varphi_0) & \mbox{in } \, \Omega . \label{old7}
\end{align}

In this case we can show the following  existence result.

\begin{theorem} \label{existtheoold}
 Let $T \in (0,\infty]$ and set either $I = [0,\infty)$, if $T = \infty$ or $I=[0,T]$, if $T < \infty$.
 Let Assumption \ref{assumptions} hold, $\bv_0 \in L^2_\sigma(\Omega)$ and $\varphi_0 \in H^1(\Omega)$
 with $|\varphi_0| \leq 1$ almost everywhere and $\int_\Omega \hspace*{-14pt}{-}\hspace*{5pt} \varphi_0 \, dx \in (-1,1)$. 
 Then there exists a weak solution $(\bv,\varphi,\mu)$ of \eqref{old1}-\eqref{old7} in the following sense:
 \begin{align*}
  & \bv \in BC_w(I;L^2_\sigma(\Omega)) \cap L^2(0,T;H_0^1(\Omega)^d) \,, \\
  & \p \in BC_w(I;H^1(\Omega)) \cap L^{\frac{4}{3}}_{\mbox{\footnotesize uloc}}(I;H^2(\Omega)) \,, \; \
        \Psi'(\p) \in L^{\frac{4}{3}}_{\mbox{\footnotesize uloc}}(I;L^2(\Omega)) \cap L^2_{\mbox{\footnotesize uloc}}(I;L^1(\Omega)) \,, \\
  & \mu \in L^2_{\mbox{\footnotesize uloc}}(I;H^1(\Omega)) 
    \; \mbox{ with } \; \nabla \mu \in L^2(0,T;L^2(\Omega)) 
 \end{align*}
 and the following equations are satisfied:
 \begin{align}  
  - \left(\rho \bv , \partial_t \bpsi \right)_{Q_T} 
  &+ \left( \operatorname{div}(\rho \bv \otimes \bv) , \bpsi \right)_{Q_T}
  + \left(2 \eta(\p) D\bv , D\bpsi \right)_{Q_T} \nonumber \\
  &= \left( \mu \nabla \varphi , \bpsi \right)_{Q_T} 
  + \left( \frac{|\bv|^2}{2} \nabla \rho , \bpsi \right)_{Q_T}  \label{weaklineold1} 
 \end{align}
 for all $\bpsi \in \left[C_0^\infty(\Omega \times (0,T))\right]^d$ with $\operatorname{div} \bpsi = 0$,
 \begin{align} 
  - \left(\p , \partial_t \zeta \right)_{Q_T} 
  + \left( \bv \cdot \nabla \p , \zeta \right)_{Q_T}
  &= - \left(m(\p) \nabla \mu , \nabla \zeta \right)_{Q_T} \,, \label{weaklineold2} 
 \end{align}
 for all $\zeta \in C_0^\infty((0,T);C^1(\overline{\Omega}))$,
 \begin{align}
   a(\varphi)^{-\frac{1}{2}} \left( \mu + \frac{\partial \rho}{\partial \varphi} \frac{|\bv|^2}{2} \right) + \widetilde{\kappa} A(\varphi)
     &= \widetilde{\Psi}_0'(\varphi)  - \Delta A(\varphi) 
     \; \mbox{ almost everywhere in } \, Q_T \; \mbox{and}  \label{weaklineold3} \\
  \left.\left( \bv,\p \right)\right|_{t=0} &= \left( \bv_0 , \p_0 \right) \,. \label{weaklineold4}
 \end{align}
 Moreover,
 \begin{align} 
  E_{\mbox{\footnotesize tot}}(\p(t),\bv(t)) &+ \int_{Q_{(s,t)}} 2 \eta(\p) \, |D\bv|^2 \, d(x,\tau) 
        + \int_{Q_{(s,t)}} m(\p) |\nabla \mu|^2 \, d(x,\tau)  \nonumber \\
   &\leq E_{\mbox{\footnotesize tot}}(\p(s),\bv(s)) \label{weaklineold5}
 \end{align}
 for all $t \in [s,\infty)$ and almost all $s \in [0,\infty)$ holds (including $s=0$).
\end{theorem}

The main ideas of the proof of Theorem \ref{existtheoold} are the same as for the previous existence result
and we will just sketch the necessary adaptations. Again we use an implicit time discretization. Let
$h = \frac{1}{N}$ for $N \in \mathbb{N}$, $\bv_k \in L^2_\sigma(\Omega)$, $\p_k \in H^1(\Omega)$ 
with $\Psi'(\p_k) \in L^2(\Omega)$ and $\rho_k = \frac{1}{2}(\tilde{\rho}_1 + \tilde{\rho}_2) 
+ \frac{1}{2} (\tilde{\rho}_2 - \tilde{\rho}_1) \p_k$ be given. We also set 
$\beta = \frac{\partial \rho}{\partial \varphi} = \frac{1}{2} (\tilde{\rho}_2 - \tilde{\rho}_1)$.
We construct $(\bv,\p,\mu)=(\bv_{k+1},\p_{k+1},\mu_{k+1})$ as solution of the 
following non-linear system.

Find $(\bv,\p,\mu)$ with $\bv \in H_0^1(\Omega)^d \cap L^2_\sigma(\Omega)$, $\varphi \in \mathcal{D}(\partial E)$
and $\mu \in H^1(\Omega)$, such that
\begin{align} 
 \left( \frac{\rho \bv - \rho_k \bv_k}{h} , \bpsi\right)_\Omega 
 &+ \left( \mbox{div}(\rho_k  \bv \otimes \bv) , \bpsi \right)_{\Omega}
 + \left(2 \eta(\p_k) D\bv , D \bpsi \right)_\Omega  \nonumber \\
 &= \left( \mu \nabla \p_k , \bpsi \right)_\Omega 
  + \left( \frac{|\bv|^2}{2} \nabla \rho_k , \bpsi \right)_\Omega  \,, \label{discreteoldline1}
\end{align}
for all $\bpsi \in C_{0,\sigma}^\infty(\Omega)$, 
%(or equivalently $\bpsi \in H^1_0(\Omega)^d \cap L^2_\sigma(\Omega)$), 
\begin{align} 
 \left( \frac{\varphi - \varphi_k}{h} , \zeta \right)_\Omega
  + \left( \bv \cdot \nabla \p_k , \zeta \right)_\Omega 
 = - \left( m(\p_k) \nabla \mu , \nabla \zeta \right)_\Omega  \label{discreteoldline2} 
\end{align}
for all $\zeta \in H^1(\Omega)$ and
\begin{align}
 & \frac{\p - \p_k}{A(\p) - A(\p_k)} \left( \mu + \beta \frac{|\bv|^2}{2} \right) 
 + \widetilde{\kappa} \, \frac{A(\p) + A(\p_k)}{2} 
 = \widetilde{\Psi}_0'(A(\p)) - \Delta A(\p) \; \mbox{ a.e. in $\Omega$.} \label{discreteoldline3}
\end{align}

To show existence of solutions of the time-discrete problem we use again the Leray-Schauder principle
and define operators $\La_{k},\F_k : X \to Y$, where
\begin{align*}
 X &= \left(H^1_0(\Omega)^d \cap L^2_\sigma(\Omega) \right) \times \D(\partial E) \times H^1(\Omega) \,, \\
 Y &= \left(H^1_0(\Omega)^d \cap L^2_\sigma(\Omega) \right)' \times \left(H^1(\Omega)\right)' \times L^2(\Omega) \,.
\end{align*}
For $\bw = (\bv,\p,\mu) \in X$ we set
\begin{align*}
 \La_k (\bw) = \begin{pmatrix}
              L_k(\bv) \\
              -\operatorname{div}(m(\p_k) \nabla \mu) + \int_\Omega \mu \, dx\\
              A(\p) + \partial \widetilde{E}(A(\p))
             \end{pmatrix} ,
\end{align*}
where 
\begin{align*}
 \left< L_k(\bv),\bpsi \right> &= \int_\Omega 2 \eta(\p_k) D\bv : D\bpsi \, dx \quad \mbox{for} \;
     \bpsi \in H^1_0(\Omega)^d \cap L^2_\sigma(\Omega) \; \mbox{ and} \\
 \left<-\operatorname{div}(m(\p_k) \nabla \mu) + \int_\Omega \mu \, dx , \eta \right> 
    &= \int_\Omega m(\p_k) \nabla \mu \cdot \nabla \eta \, dx + \int_\Omega \mu \, dx \cdot \int_\Omega \eta \, dx \quad \mbox{for} \;
     \eta \in H^1(\Omega) \,. %\\
% \left< A(\p) + \partial \widetilde{E}(A(\p)) , \eta \right> 
%    &= \int_\Omega A(\p) \, \eta 
%      + \int_\Omega \nabla A(\p) \cdot \nabla \eta
%      + \int_\Omega \widetilde{\Psi}_0'(A(\p)) \, \eta \quad \mbox{for} \; \eta \in H^1(\Omega) \,.
\end{align*}
Note that due to $\p \in \D(\partial E)$ it holds $A(\p) \in \D(\partial \widetilde{E})$ and therefore the last line 
in $\La_k(\bw)$ lies in $L^2(\Omega)$.
Furthermore for $\bw = (\bv,\p,\mu) \in X$ we define 
\begin{align*}
 \F_k(\bw) = \begin{pmatrix}
            - \frac{\rho \bv - \rho_k \bv_k}{h} - \operatorname{div}(\rho_k \bv \otimes \bv) + \mu \nabla \varphi_k + \frac{|\bv|^2}{2} \nabla \rho_k \\
            -\frac{\p - \p_k}{h} - \bv \cdot \nabla \varphi_k + \int_\Omega \mu \, dx \\
            A(\p) + \frac{\p - \p_k}{A(\p)-A(\p_k)} \left( \mu + \beta \frac{|\bv|^2}{2} \right) + \tilde{\kappa} \frac{A(\p)+A(\p_k)}{2}
           \end{pmatrix} .
\end{align*}
Then $\bw = (\bv,\p,\mu) \in X$ is a weak solution of the time discrete problem 
\eqref{discreteoldline1}-\eqref{discreteoldline3} if and only if $\La_k (\bw) - \F_k(\bw) = 0$. 
To get a continuous and even compact operator, we introduce for $0 < s < \frac{1}{4}$ the following Banach spaces
\begin{align*}
 \widetilde{X} &\ZuWeis \left( H_0^1(\Omega)^d \times L_\sigma^2(\Omega) \right) 
                  \times H^{2 - s}(\Omega) \times H^1(\Omega) \,, \\
 \widetilde{Y} &\ZuWeis L^{\frac{3}{2}}(\Omega)^d \times L^2(\Omega) \times W^{1}_{\frac{3}{2}}(\Omega) \,.
\end{align*}
With analogue arguments as above we can show that the inverse $\La_k^{-1} : \widetilde{Y} \to \widetilde{X}$ is 
a compact operator and that $\F_k : \widetilde{X} \to \widetilde{Y}$ is 
continuous and maps bounded sets into bounded sets. 

For these operators it is possible to verify the assumptions of the Leray-Schauder principle and therefore get 
a weak solution to the time-discrete problem, which additionally fulfills the energy estimate \eqref{discreteenergyestimate}.

For the remaining part about compactness in time we can show the same statements as in Lemma \ref{derivativeforchempot}
and derive then together with the fact $\int_\Omega \beta \frac{|\bv|^2}{2} \, dx \in L^\infty(0,\infty)$ from the 
energy estimate the same bounds as in \eqref{timediscrbounds}. But in this case the additional term 
$|\bv|^2$ in the line \eqref{timeintegratedline3} for the chemical potential is just bounded in the 
space $L^{\frac{4}{3}}(J;L^2(\Omega))$. By applying results of Abels and Wilke \cite{AW07} we get here
$\varphi \in L^{\frac{4}{3}}_{\mbox{\footnotesize uloc}}(J;H^2(\Omega))$ and 
$\Psi'(\p) \in L^{\frac{4}{3}}_{\mbox{\footnotesize uloc}}(J;L^2(\Omega))$. The strong convergence
\eqref{strongphiinL2H1} can then be derived by using the embedding $L^{\frac{4}{3}}(0,T;H^2(\Omega))
\hookrightarrow L^1(0,T;H^2(\Omega))$ and the result from interpolation theory 
$L^1(0,T;H^2(\Omega)) \cap L^\infty(0,T;H^1(\Omega)) \hookrightarrow L^2(0,T;H^{\frac{3}{2}}(\Omega))$,
which shows that $\varphi^N$ is bounded in the latter space and 
is enough to conclude \eqref{strongphiinL2H1} with the help of the Lemma of Aubins-Lions.

To show strong convergence of $\bv^N$ to $\bv$ in $L^2(\Omega_T)$ in this case, 
we have to bound instead of the term $\bv^N \otimes \nabla \mu^N$ in the first model now the term
$|\bv^N|^2 \nabla \varphi^N_h$. From the bound of $\bv^N$ in $L^2(0,T;L^6(\Omega))$ and $L^\infty(0,T;L^2(\Omega))$ 
we get that $|\bv^N|^2$ is bounded in $L^1(0,T;L^3(\Omega)) \cap L^\infty(0,T;L^1(\Omega))$. 
With the help of the interpolation results \eqref{interpolembedding} and \eqref{realinterpollebesgue} this leads to 
a bound of $|\bv^N|^2$ in $L^{\frac{4}{3}}(0,T;L^2(\Omega))$ and
together with the bound of $\nabla \varphi^N_h \in L^\infty(0,T;L^2(\Omega))$
we get in this case the boundedness of $|\bv^N|^2 \nabla \varphi^N_h$ in $L^{\frac{4}{3}}(0,T;L^1(\Omega))$.
Anyhow, we have a time integrability greater than $1$, which is enough to finish the proof.

\section*{Acknowledgement} 
This work was supported by the SPP 1506 "Transport Processes
at Fluidic Interfaces" of the German Science Foundation (DFG) through
the grant GA 695/6-1. The support is gratefully acknowledged.


\begin{thebibliography}{breitesteee}
        \bibitem[Abe07]{Abe07}Abels H., \textit{Diffuse Interface Models for Two-Phase Flows of Viscous 
          Incompressible Fluids}, Habilitation thesis, Leipzig 2007.
        \bibitem[Abe09a]{Abe09a}Abels H., \textit{Existence of weak solutions for a diffuse interface model for viscous, incompressible
          fluids with general densities}, Comm. Math. Phys., vol. 289 (2009), p.45-73.
        \bibitem[Abe09b]{Abe09b}Abels H., \textit{On a diffuse interface model for two-phase flows of viscous, incompressible fluids with
          matched densities}, Arch. Rat. Mech. Anal., vol. 194 (2009), p.463-506.
        \bibitem[Abe11]{LTModelShortTime}Abels H., \textit{Strong well-posedness of a  diffuse interface model for a viscous, 
          quasi-incompressible two-phase flow}, SIAM J. Math. Anal., vol. 44(1) (2012), p.316-340.
        \bibitem[AGG10]{AGG10}Abels H., Garcke H., Gr\"un G., \textit{Thermodynamically consistent diffuse interface 
         models for incompressible two-phase flows with different densities}, preprint Nr. 20/2010 University Regensburg, (2010).
        \bibitem[AGG11]{AGG11}Abels H., Garcke H., Gr\"un G., \textit{Thermodynamically consistent, frame indifferent diffuse interface 
         models for incompressible two-phase flows with different densities}, Math. Models Meth. Appl. Sci., vol 22(3) (2011).
        \bibitem[AR09]{AR09}Abels H., R\"oger M., \textit{Existence of weak solutions for a non-classical sharp interface
         model for a two-phase flow of viscous, incompressible fluids}, Ann. Inst. H. Poincar\'{e} Anal. Non Lin\'{e}aire, vol. 26(6) (2009),
         p.2403-2424.
        \bibitem[AW07]{AW07}Abels H., Wilke M., \textit{Convergence to equilibrium for the Cahn-Hilliard equation with a 
          logarithmic free energy}, Nonlin. Anal., vol. 67 (2007), p.3176-3193.
        \bibitem[Ama95]{Ama95}Amann H., \textit{Linear and Quasilinear Parabolic Problems, Vol. 1: Abstract Linear Theory}, 
          Birkh\"auser, Basel, 1995.
        \bibitem[AMW98]{AMW98}Anderson, D.-M., McFadden, G.B., Wheeler, A.A., \textit{Diffuse interface methods in fluid mechanics}, 
         Annu. Rev. Fluid Mech., vol. 30. Annual Reviews, Paolo Alto, CA, 1998, p.139-165.
        \bibitem[BL76]{BL76}Bergh J., L\"ofstr\"om J., \textit{Interpolation Spaces}, Springer, Berlin-Heidelberg-New York, 1976.
        \bibitem[Boy99]{Boy99} %%% Neu
          Boyer~F., \textit{Mathematical study of multi-phase flow under shear through order
          parameter formulation}, { Asymptot. Anal.}, 20(2) (1999), p 175-212.
        \bibitem[Boy01]{Boy01}Boyer F., \textit{Nonhomogeneous Cahn-Hilliard fluids},
           Ann. Inst. H. Poincaré Anal. Non Linéaire, vol. 18(2) (2001), p 225-259.
        \bibitem[CH58]{CahnHilliard}
           Cahn, J.W., Hilliard, J.E., 
           \textit{Free energy of a nonuniform system. {I}. {I}nterfacial energy},{ J. Chem. Phys.}, 28, No. 2 (1958), p. 258-267.
        \bibitem[DU77]{DU77}Diestel J., Uhl Jr. J.J., \textit{Vector Measures}, Amer. Math. Soc., Providence, RI, 1977.
        \bibitem[DSS07]{DSS07}Ding H., Spelt P.D.M., Shu C., \textit{Diffuse interface model for incompressible two-phase flows with 
           large density ratios}, J. Comp. Phys., vol. 22 (2007), p 2078-2095. 
        \bibitem[GPV96]{GPV96}Gurtin M.E., Polignone D., Vi\~{n}als J., \textit{Two-phase binary fluids and immiscible fluids described
           by an order parameter}, Math. Models Meth. Appl. Sci., vol. 6(6)  (1996), p.815-831. 
        \bibitem[HH77]{HH77}Hohenberg P.C., Halperin B.I., \textit{Theory of dynamic critical phenomena},
           Rev. Mod. Phys., vol. 49 (1977), p.435-479. 
        \bibitem[Lio69]{Lio69}Lions J.-L., \textit{Quelques M\'{e}thodes de R\'{e}solution des Probl\`{e}mes aux Limites
          Non lin\'{e}aires}, Dunod, Paris, 1969.
        \bibitem[Lio96]{Lio96}Lions P.-L., \textit{Mathematical Topics in Fluid Mechanics, Vol. 1, Incompressible Models},
          Clarendon Press, Oxford, 1996.
        \bibitem[LS03]{LS03}
           Liu C. and Shen J., \textit{A phase field model for the mixture of two incompressible fluids and
           its approximation by a {F}ourier-spectral method}, { Phys. D}, 179(3-4) (2003), p.~211-228.
        \bibitem[LT98]{LT98}Lowengrub J., Truskinovsky L., \textit{Quasi-incompressible Cahn-Hilliard fluids and topological 
           transitions}, R. Soc. Lond. Proc. Ser. A Math. Phys. Eng. Sci., vol.~454 (1998), p.2617-2654.
        \bibitem[Rou90]{Rou90}Roub\'{i}\v{c}ek T., \textit{A generalization of the Lions-Temam compact embedding theorem}, 
          \v{C}asopis P\v{e}st. Mat., vol. 115(4) (1990), p.338-342.
        \bibitem[Sho97]{Sho97}Showalter R. E., \textit{Monotone Operators in Banach Space and Nonlinear Partial Differential Equations}, 
          AMS, Providence, R.I., 1997.
        \bibitem[Sim87]{Sim87}Simon J., \textit{Compact sets in the space $L^p(0,T;B)$}, 
          Ann. Mat. Pura Appl. (4), vol.~146 (1987), p.65-96.
        \bibitem[Soh01]{Soh01}
          Sohr H., \textit{The {N}avier-{S}tokes Equations}, Birkh\"auser Advanced Texts: Basler Lehrb\"ucher, Birkh\"auser Verlag, Basel, 2001.
        \bibitem[Sta97]{Sta97}Starovo\u{\i}tov V.N., \textit{On the motion of a two-component fluid in the presence of 
          capillary forces}, Mat. Zametki, vol. 62(2) (1997), p.293-305, transl. in Math. Notes, vol. 62(1-2) (1997), p.244-254.
        \bibitem[St70]{St70}Stein E.M., \textit{Singular Integrals and Differentiability Properties of Functions}, 
          Princeton Univ. Press, Princeton, NJ, 1970.
        %\bibitem[SS92]{SS92}Simader C.G., Sohr H., \textit{A new approach to the Helmholtz decomposition and the Neumann problem in 
        %  $L^q$-spaces for bounded and exterior domains}, In: Mathematical problems relating to the Navier-Stokes equation, Vol. 11 of Ser. Adv.
        %  Math. Appl. Sci., pp.1-35, World Scientific Publishing, River Edge, 1992.
        \bibitem[Tri78]{Tri78}Triebel H., \textit{Interpolation Theory, Function Spaces, Differential Operators},
          North-Holland Publishing Company, Amster\-dam-New York-Oxford, 1978.
        \bibitem[Zei92]{Zei92}Zeidler E., \textit{Nonlinear Functional Analysis and its Applications I}, Springer, New York, 1992.
\end{thebibliography}
\end{document}